%% file: cresson-morin-mould.tex
\documentclass[reqno,11pt,french,english]{smfart}
\usepackage[latin1]{inputenc}
\usepackage[T1]{fontenc}
\usepackage{vmargin}  
\usepackage[french,english]{babel}
\usepackage{amsmath}
\usepackage{amsfonts}
\usepackage[mathscr]{eucal}
\usepackage{amssymb}
\usepackage{bm}  
\usepackage{amsthm}
\usepackage{smfthm}
\usepackage{amscd} 

\DeclareMathOperator{\id}{Id} 
 \DeclareMathOperator{\Exp}{Exp} 
\DeclareMathOperator{\Log}{Log}

\DeclareMathOperator{\ad}{ad}

\renewcommand{\tilde}{\widetilde}
\renewcommand{\hat}{\widehat}
\renewcommand{\phi}{\varphi}
\renewcommand{\epsilon}{\varepsilon}
\renewcommand{\leq}{\leqslant}
\renewcommand{\geq}{\geqslant}

\newcommand{\R}{\mathbb{R}} 
\newcommand{\C}{\mathbb{C}} 
\newcommand{\N}{\mathbb{N}} 
\newcommand{\T}{\mathbb{T}} 
\newcommand{\Z}{\mathbb{Z}} 

\newcommand{\A}{\mathfrak{A}}
\newcommand{\D}{\mathfrak{D}}

\newcommand{\E}{E}
\newcommand{\K}{\mathbb{K}}
\newcommand{\cal}[1]{\mathcal{#1}}
\newcommand{\inverse}[1]{\left(#1\right)^{-1}} 
\newcommand{\der}{\mathcal{D}} 
\newcommand{\ie}{\emph{i.e.\ }} 
\newcommand{\lam}{\lambda}

\newcommand{\freq}{\omega} 
\newcommand{\hot}{\text{ h.o.t. }}

\newcommand{\moule}[1]{#1^{\bullet}}  
\newcommand{\comoule}[1]{#1_\bullet}  
\newcommand{\moulexp}[1]{\mathrm{Exp}\moule{#1}}  
\newcommand{\moulefois}[2]{\left[\moule{#1}\right]_{\left(\times#2\right)}}  
\newcommand{\alphabet}[1]{#1^{*}}  

\newcommand{\Xsam}{X_{\mathrm{sam}}}
\newcommand{\Xtram}{X_{\mathrm{tram}}}
\newcommand{\Xlin}{X_{\mathrm{lin}}}
\newcommand{\Xnor}{X_{\mathrm{nor}}}
\newcommand{\Xres}{X_{\mathrm{res}}}
\newcommand{\Pran}{\mathrm{Pran}}
\newcommand{\Sam}{\mathrm{Sam}}
\newcommand{\Tram}{\mathrm{Tram}}

\renewcommand{\i}{\mathrm{i}\,} 
\newcommand{\longueur}[1]{\ell\left(#1\right)} 
\newcommand{\simp}{\mathrm{simp}} 
\newcommand{\norme}[1]{\lVert#1\rVert} 
\newcommand{\vabs}[1]{\lvert#1\rvert}  
\newcommand{\inv}[1]{#1^{-1}}  
\newcommand{\mot}[1]{\bm{#1}}   

\newcommand{\dep}{\partial}  

\newcommand{\eqn}[1]{\begin{equation}#1\end{equation}} 
\newcommand{\eq}[1]{\begin{equation*}#1\end{equation*}} 

\newcommand{\expe}[1]{\mathrm{e}^{\mathrm{i}#1}}  
\newcommand{\poisson}[2]{\left\{#1,#2\right\}}  

\newcommand{\algebreformelle}[2]{#1\bigl\langle\bigl\langle#2\bigr\rangle\bigr\rangle}

\newcommand{\Res}{\mathrm{Res}}

\numberwithin{equation}{section}

\theoremstyle{remark}

\renewcommand{\theenumi}{(\alph{enumi})}

\theoremstyle{theo}
\newtheorem{notation}{Notation}

\everymath{\displaystyle}

\title{Mould Calculus for Hamiltonian vector fields}
\author{Jacky Cresson}
\address{LMA\\ Université de Pau et des pays de l'Adour\\ avenue de l'Université
BP 1155\\ 64013 PAU Cedex}
\address{IMCCE\\ 77 avenue Denfert-Rochereau\\ 75014 Paris}
\email{jacky.cresson@univ-pau.fr}
\urladdr{http://www.univ-pau.fr/~jcresson/}
\author{Guillaume Morin}
\address{IMCCE\\ 77 avenue Denfert-Rochereau\\ 75014 Paris}
\address{CEREMADE\\ Université Paris-Dauphine\\ Place du Maréchal De Lattre de Tassigny\\ 75116 Paris}
\email{morin@ceremade.dauphine.fr}
\urladdr{http://www.ceremade.dauphine.fr/~morin/}

\font\xrm=wncyr10
\def\shu{\hbox{\xrm sh} }

\begin{document}
 \frontmatter
 %
 %
  \begin{abstract}
   \input{abstract}
  \end{abstract}
%
  
\maketitle
\tableofcontents
\mainmatter
\section{Introduction}

\input{intro}
\section{Reminder about moulds}
%
\input{reminder}
\section{Continuous prenormal forms of a vector field}
%
\input{continuous}
\section{Effective aspects of continuous prenormal forms}
%
\input{effective}
\section{A first approach to the Poincaré-Dulac normal form}
%
\input{normalform}
\section{The trimmed form}
%
\input{trimmed}
\section{The Hamiltonian case}
\input{hamiltonian}
\section{Kolmogorov Theorem}
%
\input{kolmogorov}

\section{Conclusion}
%
\input{conclusion}
\backmatter
\bibliographystyle{smfplain}
\bibliography{article_moules}
\end{document}

%% file: abstract.tex
 We present the general framework of Écalle's moulds in the case of linearization of a formal vector field without and within resonances. We enlighten the power of moulds by their universality, and calculability. We modify then Écalle's technique to fit in the seek of a formal normal form of a Hamiltonian vector field in cartesian coordinates. We prove that mould calculus can also produce successive canonical transformations to bring a Hamiltonian vector field into a normal form. 
 
 We then prove a Kolmogorov theorem on Hamiltonian vector fields near a diophantine torus in action-angle coordinates using moulds techniques.

%% file: intro.tex
  We deal in this text with formal normal forms for formal vector fields on $\C^\nu$. We use the \emph{mould formalism} by Jean Écalle to obtain those. The idea in this formalism is to consider vector fields as \emph{derivations} on the algebra of formal power series $\C[[x]]$ and work in the general free Lie algebras framework associated to the algebra build on these derivations. It was developed by Écalle (see~\cite{ecalle:1,ecalle:2,ecalle:3}) but didn't get the success it deserved yet. This text comes back on Écalle's idea with some precise calculus we didn't find in his works, although it was said to be right. The Hamiltonian parts (sections~\ref{section:hamilton} and~\ref{section:kolmogorov}) were also evoked by Écalle in~\cite{ecalle:4} but still not written: we hope to give here a little contribution to his work and an educational aspect. 
  
  In order to make the reader familiar with moulds and mould calculus in the search for formal normal forms, we recall in sections~\ref{section:reminder} to~\ref{section:trimmed} some of Écalle's work, and set a global framework for moulds, which is the general free Lie algebras framework. Then, in sections~\ref{section:hamilton} and~\ref{section:kolmogorov} we present already known results, with the new techniques of moulds.
  
The search for formal normal forms of vector fields has a great first theorem from a great mathematician: the linearization theorem by Poincaré. We give here a "moulds proof" of this theorem, which obviously make the small divisors appear, and moreover, arouses a universal character of moulds: the linearization mould only depend on the graduation (\ie the decomposition) of the vector field $X$. This is of great interest, because when the vector field is modified, the linearization mould is stil the same, as long as the graduation of the vector field is the same.

 The plan of this text is the following: sections \ref{section:reminder} to \ref{section:trimmed} are of pedagogical interest, and summarize the main definitions, results and techniques of Écalle's moulds we need. Most of it can be found in \cite{calcul_moulien, cresson_raissy,  ecalle:1, ecalle:2, ecalle:3, ecalle:0}. The original work we did can be found in the last two sections. More precisely:

  Section~\ref{section:reminder} recalls some basic definitions and results about mould formalism. In section~\ref{section:continuous} we define the main object of our concern: a prenormal form. That is, a vector field $X$ being given with a fixed diagonal linear part $\Xlin$, we look for a change of variables which brings $X$ into $\Xnor$, such that $[\Xlin,\Xnor]=0$.
  
   Section~3 deals with \emph{continuous} prenormal forms, following Écalle's terminology; that is, how does a prenormal form $\Xnor$ behave when the vector field $X$ is modified, its linear part being untouched? We give here a first application of the power of the mould formalism, calculating a direct transform of linearization of $X$, according to Poincaré's linearization theorem.
   
   The case of resonant vector fields rises in the next section~4: we obtain an analogous result of the classical Poincaré-Dulac theorem; nevertheless the prenormal form calculated here is not the Poincaré-Dulac normal form; Écalle calls it the \emph{trimmed} form.
   
   The last two sections focus on Hamiltonian vector fields, which was the original goal of this text: we make here a slightly modification in Écalle's formalism: where homogeneous differential operators were used, we need another graduation (\ie decomposition) of the vector field $X$ to prove that it is possible to make successive \emph{canonical} transformations to bring a formal Hamiltonian with a resonant linear part in cartesian coordinates into a trimmed form, preserving the Hamiltonian character at each step.
    
    Then, in section~7, following~\cite{giorgilli} we prove a formal Kolmogorov theorem on a formal Hamiltonian near a diophantine torus using techniques shown in section~\ref{section:trimmed}. We study here perturbations in action-angle coordinates as trigonometric polynomials for technical reasons.

%% file: reminder.tex
%
\label{section:reminder}
 All proofs and details about this section
can be found in~\cite{cresson_raissy}.
 We denote by $A$ an alphabet, finite or not, which is a semigroup for a law $\star$. In this section, a
letter of $A$ is denoted by $a$.
 $A^*$ denotes the set of all \emph{words} $\mot{a}$ build on $A$ \ie the totally ordered
 sequences $a_1\dotsb a_r, r\geq0$,
   with $a_i$ in $A$ and $r=\longueur{\mot{a}}$ the length of the word $\mot{a}$. We set the
   convention that a word of length $0$ is the empty word $\emptyset$. Moreover, $A^*_r$ denotes the set of words of exact length $r$.

  The natural operation on $A^*$ is the usual \emph{concatenation} of two words $\mot{a}$
   and $\mot{b}$ of
   $A^*$, which glues the word $\mot{a}$ to the word $\mot{b}$, \ie $\mot{a}\bullet\mot{b}$,
   or often simply
   $\mot{ab}$ when there is no ambiguity. Moreover, as $A$ is a semigroup, we define
    $\norme{\mot{a}}_\star$
    as the letter $a_1\star\cdots\star a_r$ of $A$, if $r=\longueur{\mot{a}}$.
%
  %
  %
  Now here is the main "new" object we focus on:
 \begin{defi}
 Let $\K$ be a ring, or a field, and $A$ an alphabet. A $\K$-valued mould $\moule{M}$ on $A$ is a map from
  $A^*$ to $\K$; the evaluation of the mould $\moule{M}$ on a word $\mot{a}$ is denoted by $M^{\mot{a}}$.
 \end{defi}
 The first important thing, is the natural one-to-one correspondence between moulds and non-commutative
  formal power series.
\subsection{Moulds and formal power series}
For $r\geq 0$, remember that $A^*_r$ is the set of words of length
$r$, with the convention that $A^*_0 =\{ \emptyset \}$. We denote by
$\K \langle A\rangle$ the set of finite $\K$-linear combinations of
elements of $A^*$, \ie \emph{non-commutative} polynomials on $A$
with coefficients in $\K$, and by $\K_r \langle A\rangle$ the set of
$\K$-linear combinations of elements of $A^*_r$, \emph{ i.e.} the
set of non-commutative homogeneous
 polynomials of
degree $r$. We have a natural \emph{graduation} on $\K \langle
A\rangle$ by the length of words: \eq{ \K \langle A\rangle
=\bigoplus_{r=0}^{\infty} \K_r \langle A\rangle . } The completion
of $\K \langle A\rangle$ with respect to the graduation by length
denoted by $\K \langle\langle A\rangle\rangle$ is the set of non-commutative formal
power series with coefficients in $\K$. An element of $\K
\langle\langle A\rangle\rangle$ is denoted by \eq{ \sum_{\mot{a} \in
A^*} M^{\mot{a}} \mot{a} ,\ \ M^{\mot{a}} \in \K , } where this sum
must be understood as
 \eq{ \sum_{r\geq 0} \left ( \sum_{\mot{a}\in
A^*_r } M^{\mot{a}} \mot{a} \right ) .
}
 Let $\moule{M}$ be a $\K$-valued mould on $A$; its generating series denoted by $\Phi_M$
belongs to
 $\K \langle\langle A\rangle\rangle$ and is defined by
\eq{
  \Phi_M =\sum_{\mot{a} \in A^*} M^{\mot{a}} \mot{a},
}
 or in a condensed way as $\sum_{\bullet} \moule{M} \bullet $. This
correspondence provides a \emph{one-to-one mapping} from the set of
$\K$-valued moulds on $A$, denoted by ${\cal{M}}_{\K} (A)$, and $\K
\langle\langle A\rangle\rangle$.
\subsubsection{Moulds algebra}
The set of moulds ${\cal{M}}_{\K} (A)$ inherits a \emph{structure of
algebra} from
 $\K \langle\langle A\rangle\rangle$. The sum and
product of two moulds $\moule{M}$ and $\moule{N}$ are denoted by
$\moule{M} +\moule{N}$ and $\moule{M} \cdot \moule{N}$ or
$\moule{M}\times\moule{N}$ respectively and defined by
 \eq{
  \begin{aligned}
   (\moule{M} +\moule{N} )^{\mot{a}} & =  M^{\mot{a}} +N^{\mot{a}} ,\\
   (\moule{M}\times\moule{N})^{\mot{a}} =(\moule{M} \cdot \moule{N} )^{\mot{a}}
         & =
                \sum_{\mot{a}^1 \mot{a}^2 =\mot{a}} M^{\mot{a}^1} N^{\mot{a}^2} ,
   \end{aligned}
 }
for all $\mot{a} \in A^*$ where this latter sum corresponds to all
the partitions of $\mot{a}$ in two words $\mot{a}^1$ and $\mot{a}^2$
of $A^*$. The product operation is associative.

The neutral element for the mould product is denoted by $\moule{1}$ and defined by
 \eq{
  \moule{1} = \begin{cases}
   1 & \text{ if $\bullet=\emptyset$,}\\
   0 & \text{ otherwise.}\end{cases}
  }
Let $\moule{M}$ be a mould. We denote by $\inverse{\moule{M}}$ the inverse of
$\moule{M}$ for the mould product when it exists, \ie the solution of the mould
equation:
\eq{
\moule{M} \cdot \inverse{\moule{M}} =\inverse{\moule{M}} \cdot \moule{M} =\moule{1} .
}
The inverse of a mould $\moule{M}$ exists if and only if $M^{\emptyset}\neq0$.
\subsubsection{Composition of moulds}
Assuming that $A$ possesses a \emph{semi-group} structure, we can
define a non-commutative version of the classical operation of
\emph{substitution} of formal power series.

We denote by $\star$ an internal law on $A$, such that $(A ,\star )$ is a semi-group. We denote by
$\norme{\,}_{\star}$ the mapping from $A^*$ to $A$ defined by
\eq{
\begin{aligned}
\norme{\,}_{\star} :   A^* & \longrightarrow   A ,\\
  \mot{a} =a_1 \dots a_r & \longmapsto   a_1 \star\dots \star a_r .
\end{aligned}
}
The $\star$ will be omitted when clear from the context.

The set $\K\langle\langle A\rangle\rangle$ is graded by
$\norme{\,}_{\star}$. A \emph{homogeneous component} of degree
$a^\prime$ of $A$, of a non-commutative serie $\Phi_{M}
=\sum_{\mot{a} \in A^*} M^{\mot{a}} \mot{a}$ is the quantity \eq{
\Phi_M^{a^\prime} =\sum_{\substack{\mot{a} \in A^*\\
\norme{\mot{a}}_{\star} =a^\prime}} M^{\mot{a}} \mot{a} . } We have
by definition \eq{ \Phi_M =\sum_{a\in A} \Phi_M^a . }
\begin{defi}[Composition]
Let $(A ,\star )$ be a semi-group structure. Let $\moule{M}$ and $\moule{N}$ be
two moulds on ${\cal{M}}_{\K} (A)$ and $\Phi_M$, $\Phi_N$ their associated
generating series. The substitution of $\Phi_N$ in $\Phi_M$, denoted by $\Phi_M \circ \Phi_N$
is defined by
\begin{equation}
\label{eq:substi2} \Phi_M \circ \Phi_N =\sum_{\mot{a} \in A^*}
M^{\mot{a}} \Phi_N^{\mot{a}} ,
\end{equation}
where $\Phi_N^{\mot{a}}$ is given by $\Phi_N^{a_1} \dots \Phi_N^{a_r}$ for $\mot{a} =a_1 \dots a_r$.\\
We denote by $\moule{M} \circ \moule{N}$ the mould of ${\cal{M}}_{\K} (A)$ such that
\begin{equation}
\label{eq:substi} \Phi_M \circ \Phi_N =\sum_{\mot{a} \in A^*}
(\moule{M} \circ \moule{N} )^{\mot{a}} \mot{a} .
\end{equation}
\end{defi}

Equation~\eqref{eq:substi} defines a natural operation on moulds denoted by $\circ$
and called \emph{composition}. Using $\norme{\,}_{\star}$ we can
give a closed formula for the composition of two moulds.
\begin{lemm}
Let $(A,\star )$ be a semi-group and $\moule{M}$, $\moule{N}$ be two moulds of
${\cal{M}}_{\K}(A)$.

For the empty word, $(\moule{M}\circ\moule{N})^\emptyset=M^\emptyset$, and for all $\mot{a} \in A^*$ of length at least $1$:
\begin{equation}
(\moule{M} \circ \moule{N} )^{\mot{a}} =\sum_{k=1}^{\longueur{\mot{a} }}
  \sum_{\mot{a}^1 \dotsc \mot{a}^k \overset{*}{=} \mot{a} } M^{\norme{\mot{a}^1}_{\star}
\cdots \norme{ \mot{a}^k}_{\star}} N^{\mot{a}^1} \cdots N^{\mot{a}^k} ,
\end{equation}
where $\mot{a}^1 \cdots \mot{a}^k \overset{*}{=} \mot{a}$ denotes all the partitions of
$\mot{a}$ such that $\mot{a}^i \neq \emptyset$, $i=1,\dots ,k$.
\end{lemm}
%

\begin{defi}\label{defi:moule_I}
The neutral element for the mould composition is denoted by
$\moule{I}$ and defined by:
 \eq{ \moule{I} =
  \begin{cases}
   1 &\text{ if $\longueur{\bullet}=1$},\\
   0 &\text{ otherwise},
   \end{cases}
}
 where $\longueur{\bullet}$ denotes the length of a word of $A^*$.
\end{defi}
\subsubsection{Exponential and logarithm of moulds}
We denote by $(\K \langle\langle A\rangle\rangle )_*$ the set of non\--commutative
formal power series without a constant term. We define the
\emph{exponential} of an element $x\in \K \langle\langle
A\rangle\rangle$ , denoted by $\exp (x)$ using the classical
formula:
\eq{
\exp (x)=\sum_{n\geq 0} \frac{x^n}{n!} .
 }
 The \emph{logarithm} of an element $1+x\in 1+(\K \langle\langle
A\rangle\rangle )_*$ is denoted by $\log (1+x)$ and defined by
\eq{\log (1+x )=\sum_{n\geq 0} (-1)^{n+1} \frac{x^n}{n} .}
 These two applications have their natural counterpart in ${\cal{M}}_{\K}
(A)$.
\begin{defi}
Let $\moule{M}$ be a mould of ${\cal{M}}_{\K} (A)$ and $\Phi_M$ the
associated generating series. Assume that $\exp (\Phi_M )$ is
defined. We denote by $\Exp \moule{M}$ the mould satisfying the
equality
\eq{ \exp \left( \sum_{\bullet} \moule{M} \bullet \right )
=\sum_{\bullet} \Exp \moule{M} \bullet .}
\end{defi}
Simple computations lead to the following direct definition of
$\Exp$ on moulds:
  \eq{
  \Exp \moule{M} =\sum_{n\geq 0}
  \frac{\moulefois{M}{n}}{n!} ,
  }
 where $\moulefois{M}{n}$, for $n$ in $\N^*$, stands for
 \eq{
   \moulefois{M}{n} = \underbrace{\moule{M}\times \dotsm
   \times\moule{M}}_{n\text{ times}} 
}
and $\moulefois{M}{0} = \moule{1}$.
The same procedure can be applied to define the logarithm of a mould.
\begin{defi}
Let $\moule{M}$ be a mould of ${\cal{M}}_{\K} (A)$ and $\Phi_M$ the
associated generating series. Assume that $\log (1+\Phi_M )$ is
defined. We denote by $\Log \moule{M}$ the mould satisfying the
equality
\eq{ \log \left( 1+ \sum_{\bullet} \moule{M} \bullet \right
) =\sum_{\bullet} \Log \moule{M} \bullet . }
\end{defi}
A direct definition of $\Log$ is then given by
\eq{
 \Log \moule{M}
=\sum_{n\geq 0} (-1)^{n+1} \frac{\moulefois{M}{n}}{n!} .
 }
As $\exp$ and $\log$ satisfy $\exp \circ \log =\log \circ \exp =\id$,
we have \eq{ \Exp \left ( \Log \moule{M} \right ) =\Log \left ( \Exp
\moule{M} \right ) =\moule{1} . }

\subsection{Algebraic aspects of moulds}
\input{rajout_2}

%% file: rajout_2.tex
Let $A$ be a set, we recall the concept of \emph{free Lie algebra} denoted by $\cal{L}_{\K} (A)$ generated by $A$, see \cite{reu},\cite{serre}. Its
elements are formal expressions in Lie brackets $[. ,. ]$ of $A$ subject only to the Jacobi identity. We have $\cal{L}_{\K} (A)\subset\algebreformelle{\K}{A}$, the \emph{enveloping algebra} of $\cal{L}_{\K} (A)$. We denote by $\cal{G}_{\K} (A)$ the set of automorphisms of
$A$, \emph{i.e.} the \emph{Lie group} associated to $\cal{L}_\K (A)$.\\

Let $\Delta : \algebreformelle{\K}{A}\rightarrow \algebreformelle{\K}{A} \otimes \algebreformelle{\K}{A}$ be the algebra morphism defined for all $a\in A$ by $\Delta (a)=a\otimes 1 + 1\otimes a$ and extended to $\algebreformelle{\K}{A}$ by linearity. Using $\Delta$ we can
characterize the element of $\algebreformelle{\K}{A}$ belonging to $\cal{L}_{\K} (A)$.
\begin{defi}
An element $P \in \algebreformelle{\K}{A}$ is called primitive if $\Delta (P)=P\otimes 1 +1\otimes P$, and group-like if $\Delta (P)=P\otimes P$.
\end{defi}
An important result on free Lie algebras is:

\begin{lemm}
\label{caracterisation-primitif}
An element $P\in \cal{L}_{\K} (A)$ if and only if $P$ is primitive and $P\in \cal{G}_{\K} (A)$ if and only if $P$ is group-like.
\end{lemm}

This result can be formulated directly on coefficients leading to two symmetries for moulds on $A$. In order to state the result, we define the combinatorial notion of \emph{shuffle} product on $A^*$:

\begin{defi}
The shuffle product denoted by $\shu :A^* \times A^* \rightarrow P(A^* )$ is defined inductively on $A^*$ by $\mot{a}\shu e =e\shu \mot{a}=\mot{a}$ and $x\mot{a} \shu y\mot{b} =
x( \mot{a} \shu y\mot{b} ) \bigcup y (x\mot{a} \shu \mot{b} )$ for all $x,y\in A$, $\mot{a},\mot{b} \in A^*$.
\end{defi}

As an example, if $\mot{a}=(x_1,x_2 )$ and $\mot{b}=(x_3)$, we have
\eq{\mot{a}\shu \mot{b}=\{ (x_1,x_2,x_3),(x_1,x_3,x_2),(x_3,x_1,x_2) \}.}

According to~\cite{reu}, we have
\eq{\Delta \left ( \sum_{\mot{a} \in A^*} { M}^{\mot{a}} \mot{a} \right ) =\sum_{\mot{a},\mot{b}\in A^*} \left ( \displaystyle\sum_{\mot{c} \in \mot{a}
\shu \mot{b}} {M}^{\mot{c}} \right ) \mot{a} \otimes \mot{b} .}

Lemma~\ref{caracterisation-primitif} can be formulated as follow:

\begin{lemm}
An element $P\in \algebreformelle{\K}{A}$, $P=\sum_{\mot{a} \in A^*} {P}^{\mot{a}} \mot{a}$ is primitive (resp. group-like) if and only if
\eq{\sum_{\mot{c} \in \mot{a} \shu \mot{b}} { M}^{\mot{c}} =0\ \ \forall \mot{a},\mot{b} \in A^* \setminus \{ e\}, \eqno{(\star )}.}
\eq{\left (\displaystyle\sum_{\mot{c} \in \mot{a} \shu \mot{b}} { M}^{\mot{c}} ={ M}^{\mot{a}} { M}^{\mot{b}} \ \ \forall \mot{a},\mot{b} \in A^*,\ \ \
{\rm resp.} \right ) \eqno{(\star \star )}.}
\end{lemm}

We now introduce Ecalle's terminology for moulds corresponding to primitive or group-like elements in $\algebreformelle{\K}{A}$

\begin{defi}
A mould $\moule{M} \in \cal{M}_{\K} (A)$ is called alternal (resp. symetral) if $\moule{M}$ satisfies ($\star$) (resp. ($\star \star$)).
\end{defi}

A primitive element of $\algebreformelle{\K}{A}$ being given, we can easily obtain the corresponding element of $\cal{L}_{\K} (A)$. We denote by $\cal{I}$ the ideal of $\K \langle A\rangle$ generated by all polynomials without a constant term.
We denote by $\psi :\cal{I} \rightarrow \cal{L}_{\K} (A)$ the mapping defined for all $\mot{a} \in A^*$,
$\mot{a}=(a_1 ,\dots ,a_r)$ by
\eq{\psi (\mot{a} )= \frac{1}{r} [[ \dots [[a_1 ,a_2],a_3],\dots ], a_r] ,}
and extended by linearity to $\cal{I}$.

According to \cite{serre}, we have the following result called the \emph{projection lemma} by Ecalle:

\begin{lemm}\label{lem:projection}
We have $\psi\mid_{\cal{L}_{\K} (A)} =\id_{\cal{L}_\K (A)}$.
\end{lemm}

As a consequence, for an alternal mould $\moule{M} \in \cal{M}_\K (A)$, we have
\eq{\sum_{\mot{a}\in A^*} {M}^{\mot{a}} \mot{a} =\sum_{\mot{a} \in A^*} {M}^{\mot{a}} \psi (\mot{a} )=\sum_{r\geq 1}
{1\over r} \sum_{\mot{a} \in A^*_r } {M}^{\mot{a}} [[\dots [[a_1 ,a_2],a_3],\dots ], a_r].}

%% file: continuous.tex
%
\label{section:continuous}
From now on, $\nu$ will be an integer, and $X$ a vector field on
$\C^\nu$ such that $X(0)=0$. We want to obtain some particular form
of $X$ within a change of variable.
 If $m=(m_1,\dotsc,m_\nu)\in\Z^\nu$, we use the notation $x^m$ for $x_1^{m_1}\cdots x_\nu^{m_\nu}$ and
 $\dep_{x_i}$ for $\frac{\dep}{\dep x_i}$.
  \begin{defi}
  A \emph{differential operator} is an element of $\C[[x]] [[\dep_{x_1},\dotsc,\dep_{x_\nu}]]$
   \ie a formal power series in the $\dep_{x_i}$ whose coefficients are (commutative) formal power series in $x$.

   The \emph{order} of a differential operator is the degree of the corresponding polynomial
   in the variables $\dep_{x_1},\dotsc,\dep_{x_\nu}$.

  If $n$ is in $\Z^\nu$, a \emph{homogeneous differential operator of degree $n$}
   is a differential operator $B_n$ such that for all $m$
    in $\N^\nu$, there exists a $\beta_{n,m}$ in $\C$, such that:
  \eq{ B_n(x^m) = \beta_{n,m}x^{n+m}.}
  \end{defi}
We usually omit the composition operator $\circ$ when composing
homogeneous differential operators: we write $B_{n^1}\dotsm B_{n^r}$
for $B_{n^1}\circ\dotsm\circ B_{n^r}$. Moreover, we denote
$B_{\mot{n}}$ for $B_{n^1}\dotsm B_{n^r}$ where $\mot{n}$ is the
word $\mot{n}=n^1\dotsb n^r$.

\begin{rema}\label{rema:diff}  Finally, remark that if $B_n$ is a differential homogeneous
  operator of order $1$ and degree $n$, $B_n = \sum_{i=1}^\nu B_n(x_i)\dep_{x_i}$.
 When there is no ambiguity, we denote $\sum_\bullet\moule{M}\comoule{B}$ for
 $\sum_{\mot{n}\in A(X)^*} M^{\mot{n}} B_{\mot{n}}$.
\end{rema}
\subsection{Resonant normal form}\label{subsec:resonant}
 Now consider a vector field $X=\sum_{i=1}^\nu X_i(x) \dep_{x_i}$ on
 $\C^\nu$ (with $X(0)=0)$; it is always possible to write this vector field as
 \eq{ X = \Xlin + \sum_{n\in A(X)} B_n,}
 where the $B_n$ are homogeneous differential operators of degree
 $n$ and order $1$, $A(X)$ is an appropriate part of $\Z^\nu$ (that will be the
 alphabet) and $\Xlin$ the linear part.
As an example, for the following vector field on $\C^2$~:
\eqn{X=\lambda_1 x\dep_x+\lambda_2 y\dep_y+(a_{20}x^2+a_{11}xy+a_{02}y^2)\dep_x %
                       +(b_{20}x^2+b_{11}xy+b_{02}y^2)\dep_y,}
 we have $\Xlin=\lambda_1 x\dep_x+\lambda_2 y\dep_y$ and the
 homogeneous differential operators are:
\begin{equation*}
\begin{aligned}
B_{(1,0)} &= x(a_{20}x\dep_x+b_{11}y\dep_y),\\
B_{(0,1)} &= y(a_{11}x\dep_x+b_{02}y\dep_y),\\
B_{(-1,2)} &= a_{02}y^2\dep_x,\\
B_{(2,-1)} &= b_{20}x^2\dep_y.
\end{aligned}
\end{equation*}
The alphabet here is $A(X)=\{(1,0),(0,1),(-1,2),(2,-1)\}$.
   The linear part of $X$ is always supposed to be of a diagonal
form (see~\cite{martinet} for instance); we have then:
   \eq{\Xlin = \sum_{i=1}^\nu \lam_i x_i \dep_{x_i},}
where $\lam=(\lam_1,\dotsc,\lam_\nu) \in\C^\nu$ is the spectrum of
$\Xlin$.
\begin{rema}
 We use here the graduation by degree for the operators $B_n$ but it is not unique;
for instance, let us set $\Omega=\{ \lambda\cdot n, n\in A(X)\}$ and
 $\mathbb{B}_\omega = \sum_{\substack{n\in A(X)\\ \lambda\cdot n=\omega}} B_n$; we still
 have $X=\Xlin + \sum_{\omega\in\Omega}\mathbb{B}_\omega$. This graduation is used by Ecalle
 in~\cite{ecalle:0} but the operators $\mathbb{B}_\omega$ are not homogeneous. We also use
 this graduation in section~\ref{section:hamilton}.
\end{rema}
\begin{defi}
When the field $X$ is written as follows:
 \eq{X=\sum_{i=1}^\nu \lam_i
  x_i \dep_{x_i} + \sum_{n\in A(X)} B_n,}
it is said to be in \emph{prepared form.}
\end{defi}
\subsection{An algebraic point of view}
Starting from this writing, we look for a change of variables $h$ in $\C[[x]]$,
from $\C^\nu$ to $\C^\nu$, which is tangent to identity (\ie
$h(x)=x+\dotsb$), to simplify $X$. We define moreover the
\emph{substitution morphism} $\Theta_h$ as:
 \eq{
  \begin{aligned}
   \Theta_h : \C[[x]] &\to \C[[x]]\\
            \phi &\mapsto \phi\circ h
   \end{aligned}
   }
It will be denoted as $\Theta$ when no ambiguity.
Remark that $h$ is a change of variables, and is then one-to-one.
Hence $\Theta$ is an automorphism of $\C[[x]]$.

 The natural action of the vector field on formal power series $\phi$
 of $\C[[x]]$ is the derivation $\phi\mapsto X\cdot\phi$ where
 $X\cdot\phi = \sum_{j=1}^\nu X_j(x)\frac{\dep\phi}{\dep x_j}(x)$.
 Hence the change of variables $h$ must let the following diagram be commutative:
 \eq{
  \begin{CD}
  \phi @>X>>X\cdot\phi \\
  @V{h}VV @VV{h}V\\
 \phi\circ h @>>{\Xnor}> \Xnor\cdot(\phi\circ h )
 \end{CD}
 }
Hence, for all $\phi\in\C[[x]]$ we must have $\Xnor\cdot(\phi\circ h) = (X\cdot\phi)\circ h$,
that is $\Xnor(\Theta\phi)=\Theta(X\phi)$, \ie $\Xnor\Theta=\Theta X$, or $\Xnor=\Theta
X\inv{\Theta}$.
 Now, the object we are looking for is the "new" field, $\Xnor$, defined
by: \eq{\Xnor = \Theta X \inv{\Theta}.}
\begin{defi}
We say that a field $\Xnor$ is a \emph{prenormal form} of $X$ if
$\Xnor$ is conjugated to $X$ and
$\bigl[ \Xnor, \Xlin \bigr]=0$ where $[\, ,]$ are the usual Lie
brackets.
\end{defi}
We recall also the following definitions from
Arnold~\cite{arnold_1}, \S~22, p.175 and p.178:
\begin{defi}\label{defi:resonant}
The spectrum $\lam=(\lam_1,\dotsc,\lam_\nu)$ is \emph{resonant} if
there is at least one $s$ in $\{1,\dotsc,\nu\}$ such that there
exists $m$ in $\N^\nu$, $\vabs{m}\geq2$ such that:
 \eq{ \lam_s = m\cdot\lam = \sum_{i=1}^\nu m_i\lam_i.}

 Moreover, a vectorial monomial $x^m\dep_{x_s}$ is resonant if
 $\lam_s = m\cdot\lam, \vabs{m}\geq2$.
\end{defi}
 For a monomial $x^m\dep_{x_i}$, with $m$ in $\N^\nu$ and $i$ in
 $\{1,\dotsc,\nu\}$, we have
   \eq{[x^m\dep_{x_i},\Xlin] = x^m\bigl(\lam_i-\sum_{j=1}^\nu m_j\lam_j\bigr)\dep_{x_i};}
 thus a prenormal form is only made of resonant monomials, \ie
  \eq{\Xnor = \Xlin + \sum_{i=1}^\nu \sum_{m\in R_i(X)} a_m x^m\dep_{x_i}}
  with $a_m\in\C$ and $R_i(X)=\{m\in\Z^\nu-\{0\}, \vabs{m}\geq2, \lambda_i=m\cdot\lambda\}$.

\subsection{Non-unicity of prenormal forms}
For an integer $k$ greater than $2$ we denote by $E_k$ the set of
all homogeneous vector fields on $\C^\nu$ of degree $k$, that is
vector fields on $\C^\nu$ whose each component is a homogeneous
polynomial in $x_1,\dotsc,x_\nu$ of degree $k$. Now, let $E$ be
defined by:
 \eq{E=\bigoplus_{k\geq2}E_k.}
Any vector field $Y$ of $E$ writes then $Y=\sum_{k\geq2} Y_k$ where
$Y_k$ is in $E_k$.

 The part $\Xlin$ being fixed, we denote by $\ad_{\Xlin}$ the
 application defined
  by $Y\mapsto\bigl[Y,\Xlin\bigr]$. Remark that $\ad_{\Xlin}(E_k)\subset
  E_k$ for any $k\geq2$, and therefore that $\ad_{\Xlin}$ may be
  defined on $E$.

\begin{prop}\label{prop:trivial}
Let $\Xlin= \sum_{i=1}^\nu\lambda_i x_i\dep_{x_i}$ and
  $\left\{\begin{aligned}\ad_{\Xlin}&:E\to E\\
                    Y&\mapsto[Y,\Xlin]\end{aligned}\right.$.
Then, $\ker(\ad_{\Xlin})=\{0\}$ if and only if $\lambda$ is non-resonant.
\end{prop}
\begin{proof}
By linearity, it is sufficient to prove it for a homogeneous polynomial of degree $k$, and
even for a monomial $x^m\dep_{x_i}$.
From definition~\ref{defi:resonant}, if $\lambda$ is non-resonant, then
 $\bigl(\lam_i-\sum_{j=1}^\nu m_j\lam_j\bigr)\neq0$ hence $\ker\ad_{\Xlin}=\{0\}$.
  Conversely, if $\ker\ad_{\Xlin}=\{0\}$ then $\bigl(\lam_i-\sum_{j=1}^\nu m_j\lam_j\bigr)\neq0$.
Hence $\lambda$ is non-resonant.
\end{proof}
\begin{coro}\label{coro:nonresonant}
If the spectrum $\lambda$ of $\Xlin$ is non-resonant, then $\Xnor=\Xlin$.
\end{coro}
\begin{proof}
Indeed, a prenormal form is a sum of $\Xlin$ and only resonant monomials.
\end{proof}

 In the case where the spectrum $\lam$ is resonant, have the following proposition:
 \begin{prop}
 If the spectrum $\lambda$ of $\Xlin$ is resonant, a prenormal form of $\Xlin$ is not unique.
\end{prop}
\begin{proof}
 Indeed, a prenormal form is defined up to a vector field $Y$ in the kernel
 of $\ad_{\Xlin}$, which is not trivial, after proposition~\ref{prop:trivial}.
 \end{proof}

  There is thus a choice to make, which could
simplify the transformation. Baider~\cite{baider} and
Gaeta~\cite{gaeta} have two interesting approaches of that.

   We denote by $\Res(\E)$ the kernel of $\ad_{\Xlin}$. If $X$ is in $\E$, when looking
   for a prenormal form of $X$, we want to write:
    \eq{ X = \Xlin + \Xres, \text{ where $\Xres$ is in $\Res(\E)$}.}
This writing infers the direct sum decomposition:
 \eq{\E = \Res(\E)\bigoplus S,}
where $S$ is a supplementary which can be chosen in many ways. A
convenient way to chose $S$ is to provide $\E$ with a scalar product
such that
 \eq{\E=\bigoplus^\bot\E_k,}
where $\E_k$ is the homogeneous component of degree $k$ of $\E$.

\subsection{Continuous prenormal forms}
 We denote by
$\algebreformelle{\C}{\mot{B}}$ for the algebra
$\algebreformelle{\C}{\bigl(B_n\bigr)_{n\in A(X)}}$ of non
commutative formal series build on the $B_n$ operators.
\begin{prop}
 There is a one-to-one correspondence between $\algebreformelle{\C}{\mot{B}}$
and $\algebreformelle{\C}{A(X)}$ given by:
\eq{ \sum_{\mot{n}\in A(X)^*}M^{\mot{n}}B_{\mot{n}} \mapsto \sum_{\mot{n}\in A(X)^*} M^{\mot{n}}\mot{n}.}
\end{prop}

As there is also a one-to-one correspondence between $\algebreformelle{\C}{A(X)}$
and $\mathcal{M}_\C\bigl(A(X)\bigr)$, we have the following writing:
  \eq{X=\Xlin + \sum_\bullet\moule{I}\comoule{B}.}
If $k=(k_1,\dotsc,k_\nu)$ is a $\nu$-uplet, we denote by
 \eq{\omega(k)=\lam\cdot k= \sum_{i=1}^\nu \lambda_i k_i}
 where $\lam$ is
the (fixed) spectrum of $\Xlin$. Remember also that if
$\mot{a}=(a^1,\dotsc,a^r)$ is in $A^*$, we write $\norme{\mot{a}} =
\norme{\mot{a}}_+ = a^1+\dotsm+a^r$.
We have the following result:
\begin{lemm}
\label{lem:crochet}
 For any word $\mot{n}$ of $A^*$, we have
$\bigl[\Xlin,B_{\mot{n}}\bigr] =
\omega(\norme{\mot{n}})B_{\mot{n}}$.
\end{lemm}
\begin{proof}
We prove the result by induction on the length $r$ of the word.
Remember that for a word $n$ of length $r=1$, $B_n = \sum_{j=1}^\nu
B_n(x_j) \dep_{x_j}$. Hence
 \eq{
  \begin{aligned}
\bigl[ \Xlin,B_n \bigr] &= \Xlin B_n - B_n \Xlin \\
              &= \sum_{j=1}^\nu\biggl( \sum_{i=1}^\nu \lam_i x_i
\dep_{x_i}\bigl[B_n(x_j)\bigr] - B_n(x_j)\lam_j\biggr) \dep_{x_j}.
\end{aligned}
}
By definition of $B_n$, homogeneous differential operator of degree
$n$, we can write $B_n(x_j)= \beta_j x^{\check{n}_j}$ where
$\beta_j$ is a complex and $\check{n}_j =
(n_1,\dotsc,n_j+1,\dotsc,n_\nu)$. Then for any $(i,j)$ in
$\{1,\dotsc,\nu\}^2$, $x_i\dep_{x_i}\bigl[B_n(x_j)\bigr] =
\beta_j(n_i+\delta_{ij})x^{\check{n}_j}$.
 Hence
\eq{ \begin{aligned}
 \bigl[ \Xlin,B_n \bigr]
   &= \sum_{j=1}^\nu \beta_j x^{\check{n}_j} \bigl( \sum_{i=1}^\nu \lam_i(n_i+\delta_{ij})-\lam_j
       \bigr) \\
   &= \sum_{j=1}^\nu \beta_j x^{\check{n}_j} \omega(n) =
   \omega(n)B_n.
   \end{aligned}
}
 Now let be $r\geq2$ fixed; we set the assumption that for every
 word $\mot{m}$ of length less than $r-1$, then $\bigl[ \Xlin,
 B_{\mot{m}} \bigr] = \omega(\norme{\mot{m}})B_{\mot{m}}$. For a
 word $\mot{n}$ of length $r$ we write $\mot{n}=\mot{m}p$ where
 $\longueur{\mot{m}}=r-1$ and $\longueur{p}=1$. Then:
 \eq{
\begin{aligned}
 \bigl[ \Xlin, B_{\mot{n}} \bigr]
  &= \bigl[\Xlin,B_{\mot{m}}B_p\bigr]\\
  &= \Xlin B_{\mot{m}}B_p-B_{\mot{m}}B_p\Xlin \\
  &= \Xlin B_{\mot{m}}B_p-B_{\mot{m}}\Xlin B_p + B_{\mot{m}}\Xlin
     B_p - B_{\mot{m}}B_p\Xlin \\
  &= \bigl[\Xlin,B_{\mot{m}}\bigr] B_p +
     B_{\mot{m}}\bigl[\Xlin,B_p\bigr] \\
  &= \omega(\norme{\mot{m}})B_{\mot{n}} + B_{\mot{m}}\omega(p)B_p \\
  &= \omega(\norme{\mot{n}}) B_{\mot{n}}.
\end{aligned}
}
\end{proof}
\begin{notation}
For a letter $n$ in $A(X)$, $\omega(n)$ stands for $\lam\cdot n$.
This notation extends to words $\mot{n}$ in $A(X)^*$, by
$\omega(\norme{\mot{n}})$. We set the convention that
$\omega(\emptyset)=0$.
\end{notation}
 \begin{defi}
 Let $\mot{n}$ be in $A(X)^*$; $\mot{n}$ is \emph{resonant} if
 \eq{\omega\bigl(\norme{\mot{n}}\bigr)=\lam\cdot \norme{\mot{n}}=0.}
 \end{defi}
 We can now define the particular forms we are looking for:
\begin{defi}
The vector field $X$ is said to be in \emph{continuous prenormal
form with respect to the alphabet $A$} if there is a change of
variable that conjugates $X$ to $\Xnor$ so that
\eq{ \Xnor = \Xlin + \sum_{\mot{n}\in A(X)^*} \Pran^{\mot{n}}
B_{\mot{n}},
     \text{ with $\Pran^{\mot{n}}=0$ if
     $\omega(\norme{\mot{n}})\neq0$.}}
\end{defi}
\begin{rema}
Notice that this definition implies that the mould $\moule{\Pran}$ must be
an alternal mould, as $\Xnor-\Xlin$ is a vector field, hence a primitive element of $\algebreformelle{\C}{\mot{B}}$. Therefore, after lemma~\ref{lem:projection}, $\Xnor-\Xlin$
is an element of the \emph{Lie algebra} spanned by the $\{B_n, n\in A(X)\}$.
We are thus trying to write
elements of $\ker\ad_{\Xlin}$ in the Lie algebra spanned by the
$B_n$. There is nevertheless no reason why all elements of that
kernel should be writing that way. However this condition is
justified as $\Xnor$ would then be intrinsic to $X$. It is therefore Ecalle's choice
to look for elements of $\ker\ad_{\Xlin}$ in that Lie algebra.
\end{rema}
\begin{rema}
  A second remark, is that this definition of a continuous prenormal form depends on the
writing of $X$ as a decomposition in the operators $B_n$, hence on
the alphabet $A(X)$. We will see in section~\ref{section:hamilton}
that some choices of decomposition may be better than others.
\end{rema}
We have then the following result:
\begin{theo}
A continuous prenormal form is a prenormal form.
\end{theo}
\begin{proof}
 The result is obtained by applying lemma~\ref{lem:crochet}
with the definition of a resonant word.
\end{proof}

%% file: effective.tex
%
%
%
%
We are looking for $\Xnor-\Xlin$ to be in the free Lie algebra of $\algebreformelle{\C}{\mot{B}}$,
hence the automorphism $\Theta$ must be in the corresponding Lie group.
That is the reason why we work in the universal enveloping algebra
$\algebreformelle{\C}{\mot{B}}$.
  Hence: \eq{\Theta = \sum_\bullet
\moule{\Theta}\comoule{B},}
with the mould $\moule{\Theta}$ being \emph{symetral} as $\Theta$ must be
an automorphism of $\C[[x]]$, \ie a group-like element of $\mathcal{L}_{\C}(A)$.
\subsection{The conjugation equation}
 The conjugation equation also writes:
 \eq{ \inv{\Theta} \Xnor \Theta = X,}
 where $\Xnor$ is the prenormal form looked for, with
  \eq{\Xnor = \Xlin + \sum_\bullet\moule{\Pran}\comoule{B}.}
The mould expansion gives:
 \eq{
  \Xlin +
\sum_\bullet \moule{I}\comoule{B} = \biggl( \sum_\bullet
\inverse{\moule{\Theta}}\comoule{B} \biggr)\biggl( \Xlin +
\sum_\bullet\moule{\Pran}\comoule{B} \biggr) \biggl( \sum_\bullet
\moule{\Theta}\comoule{B} \biggr)}
 \text{ \ie }
\eqn{\label{eq:conjugaison_1}
  \Xlin + \sum_\bullet \moule{I}\comoule{B} =
 \biggl( \sum_\bullet \inverse{\moule{\Theta}}\comoule{B}
 \biggr) \Xlin\biggl( \sum_\bullet \moule{\Theta}\comoule{B}
 \biggr) + \sum_\bullet \biggl( \inverse{\moule{\Theta}}
 \times\moule{\Pran}\times\moule{\Theta}\biggr) \comoule{B}.
}
 As we can see on this latter equation, the quantity
 $\Xlin\comoule{B}$ must be investigated.
\begin{lemm}\label{lem:xlin_1}
 Let $\moule{M}$ be a mould in $\mathcal{M}_\C\bigl(A(X)\bigr)$.
 Then:
 \eq{
  \Xlin \Bigl(\sum_\bullet \moule{M} \comoule{B}\Bigr) =
  \sum_\bullet\nabla\moule{M}\comoule{B} + \Bigl(\sum_\bullet
  \moule{M}\comoule{B} \Bigr) \Xlin,}
  where $\nabla M^{\mot{n}} = \omega(\norme{\mot{n}}) M^{\mot{n}}$
  for all $\mot{n}$ in $A(X)^*$.
 \end{lemm}
\begin{proof}
 By linearity, it is sufficient to calculate $\Xlin B_{\mot{n}}$ for a word $\mot{n}=(n^1,\ldots,n^r)$
in $A(X)^\ast$ of length $r$.
 After lemma~\ref{lem:crochet}, $\Xlin B_{\mot{n}} = \omega(\norme{\mot{n}})B_{\mot{n}}
 + B_{\mot{n}}\Xlin$.

\end{proof}
\begin{prop}
 The conjugation equation has the following mould form:
\eqn{\label{eq:conjugaison_2}
  \moule{\Theta}\times\moule{I} = \nabla \moule{\Theta} +
  \moule{\Pran}\times\moule{\Theta}.}
\end{prop}
\begin{proof}
 Using lemma~\ref{lem:xlin_1}, equation~\eqref{eq:conjugaison_1} rewrites:
  \eq{
 \moule{I} = \inverse{\moule{\Theta}} \times \nabla\moule{\Theta} +
             \inverse{\moule{\Theta}} \times \moule{\Pran} \times
             \moule{\Theta},}
 and, after left-multiplicating by $\moule{\Theta}$:
 \eq{\moule{\Theta}\times\moule{I} = \nabla \moule{\Theta} +
  \moule{\Pran}\times\moule{\Theta}.}
\end{proof}
This equation gives us a relation between the normalisator $\Theta$
and the desired prenormal form.
\subsection{The non-resonant case}
 In the case where $\lam$ is non-resonant, we must have
$\moule{\Pran}=0$ for all $\bullet$, \ie $\Xnor=\Xlin$
(linearization of $X$), after corollary~\ref{coro:nonresonant}.
We have then to solve an induction relation on
the $\moule{\Theta}$ to prove its existence

Equation~\eqref{eq:conjugaison_2} rewrites indeed:
\eq{ \moule{\Theta}\times\moule{I}=\nabla\moule{\Theta}.} Remember
definition~\ref{defi:moule_I} of $\moule{I}$ on
page~\pageref{defi:moule_I}.

 Hence for a word of length $0$, $B_{\emptyset}=\id$ and $\Theta$ is
 tangent to identity, so $\Theta^\emptyset=1$.

 For a word $n$ of length $1$, equation~\eqref{eq:conjugaison_2}
 rewrites $1=\omega(n)\Theta^n$; $\lam$ is non-resonant therefore
 $\omega(n)\neq0$, thus $\Theta^n=\frac{1}{\omega(n)}$.

 For a word $\mot{n}=(n^1,\dotsc,n^r)$ of length $r$ at least $2$,
 equation~\eqref{eq:conjugaison_2} rewrites then:
  \eq{ \Theta^{n^1,\dotsc,n^{r-1}} =
  \omega(\norme{\mot{n}})\Theta^{\mot{n}};}
 as $\lam$ is non-resonant, we have still that
 $\omega(\norme{\mot{n}})\neq0$, hence the induction formula:
  \eq{
   \Theta^{\mot{n}} =
   \frac{\Theta^{n^1,\dotsc,n^{r-1}}}{\omega(\norme{\mot{n}})},}
 hence, by induction:
 \eq{
  \Theta^{\mot{n}} = \frac{1}{\omega_1
  (\omega_1+\omega_2)\dotsm(\omega_1+\dotsb+\omega_r)},}
  where $\omega_i$ stands for $\omega(n^i)$ for $i$ in
  $\{1,\dotsc,r\}$.
We must then be sure that $\Theta$ is an automorphism of $\C[[x]]$, \ie that
$\moule{\Theta}$ is symetral. This is indeed true, see~\cite{calcul_moulien} for a proof.
\subsection{The resonant case}
In the case where $\lam$ is resonant, we set $\Pran^{\mot{n}}=0$ if
     $\omega(\norme{\mot{n}})\neq0$, to obtain a continuous
     prenormal form of $X$. However, equations on $\Theta^{\mot{n}}$ cannot
      be solved directly this time:
we take equation~\eqref{eq:conjugaison_2} and try to solve it, by
induction on the length of words.

 Remember first that $\Theta^{\emptyset}=1$, for $\Theta$ has to be tangent to
 identity; moreover, $I^{\mot{n}}=1$ if $\longueur{\mot{n}}=1$ and
 $0$ otherwise; finally we set $\Pran^{\mot{n}}=0$ if
 $\omega(\norme{\mot{n}})\neq0$.

 For the empty word $\emptyset$, equation~\eqref{eq:conjugaison_2}
 rewrites $\Theta^{\emptyset}I^{\emptyset} =
 \nabla\Theta^{\emptyset} + \Pran^{\emptyset}\Theta^{\emptyset}$ \ie
 $\Pran^{\emptyset}=0$.

For a word $n$ of length $1$, equation~\eqref{eq:conjugaison_2}
rewrites:
 \eq{\begin{aligned}
 \Theta^nI^{\emptyset} + I^n\Theta^\emptyset &= \nabla\Theta^n +
 \Pran^n\Theta^{\emptyset} + \Pran^{\emptyset}\Theta^n \\
 \iff 1 &= \omega(n)\Theta^n + \Pran^n.
 \end{aligned}
 }
 Thus, if $\omega(n)\neq0$, $\Pran^n=0$ so we can solve this equation and
 $\Theta^n=\frac{1}{\omega(n)}$. However, if $\omega(n)=0$, it gives
 $\Pran^n=1$ but we have no information on $\Theta^n$. It is
 therefore not possible to deduce all $\Theta^{\mot{n}}$ for all
 words $\mot{n}$ in $A(X)^*$ from equation~\eqref{eq:conjugaison_2}
 and the condition we set on $\moule{\Pran}$. At this moment, there
 are two possibilities:
  \begin{itemize}
  \item either add a condition, like a derivation
 condition on moulds (see~\cite{vallet} p.25); it leads to different
 prenormal forms, depending on the additional condition. This method is the "direct" method.

  \item The other possibility is an iterative method like Poincaré-Dulac.
  This is the method we expose here.
\end{itemize}
  As done usually (see~\cite{arnold_1} for instance),
when looking for a (pre)normal form, we want to write $X-\Xlin$ as a
sum of resonant terms only. Hence we set the mould $\moule{\Pran}$
this way:
 \eq{
 \Pran^{\mot{n}} =
   0 \text{ if $\mot{n}$ is non-resonant.}}

Unfortunately, it seems too complicated to solve this equation at
once.
 Nevertheless, it is possible to do a calculable step-by-step procedure to
 obtain $\moule{\Pran}$ with the wanted properties. The step-by-step
 procedure is explained in the next section.

%% file: normalform.tex
%
%
%
\subsection{The interest of being in a Lie algebra}
 The idea of the step-by-step procedure is to kill non-resonant
 terms $B_{\mot{n}}$ of order $1$ (\ie such that $\longueur{\mot{n}}=1$) at each
 step.

 We saw that the normalizator $\Theta$ is an automorphism of
$\C[[x]]$, that is a group-like element of the free Lie algebra $\mathcal{L}_\C(A)$;
 it can therefore be written as an exponential of a primitive element of this algebra,
  \ie an exponential of a vector field $V$: we thus write
 \eq{\Theta=\exp(V), \text{ where }V =
 \sum_\bullet\moule{V}\comoule{B},}
with $\moule{V}$ being alternal.
 After the mould writing of $\Theta$, and by definition of the
 exponential of a mould, we have:
  \eq{\moule{\Theta} = \moulexp{V},}
 hence equation~\eqref{eq:conjugaison_2} rewrites then:
 \eqn{\label{eq:conjugaison_3}
  \moulexp{V}\times\moule{I} = \nabla \moulexp{V} +
  \moule{\Pran}\times\moulexp{V}.}
Still we set $\Pran^{\mot{n}}=0$ if $\omega(\norme{\mot{n}})\neq0$,
but we have the same indetermination on $\moule{V}$ as we had on
$\moule{\Theta}$: we choose here to kill only non-resonant terms of
order $1$. The exponential form of the normalisator, with the help
of the Baker-Campbell-Hausdorff formula, leads us to the following
lemma:
\begin{lemm}
Let $X= \Xlin + \sum_{n\in A(X)} B_n$ a vector field in prepared
form, with fixed diagonal linear part $\Xlin$ of spectrum $\lam$.

Choosing $\Theta=\exp(V)$, where $V$ is defined by:
 \eq{
  V^{\mot{n}} = \begin{cases}
                \frac{1}{\omega(n)} &\text{ if
                $\longueur{\mot{n}}=1$ and $\omega(n)\neq0$;}\\
                0 &\text{ otherwise;}
                \end{cases}
}
 the conjugate vector field $\Xnor=\Theta\Xlin\inv{\Theta}$ has no resonant
 terms of order $1$.
\end{lemm}
\begin{proof}
After the Baker-Campbell-Hausdorff formula (see~\cite{jacobson}) we have:
\eq{ \Xnor = \sum_{n=0}^{+\infty}\frac{(-1)^n}{n!}X^{(n)},}
where $X^{(n+1)}=[X^{(n)},V]$ and $X^{(0)}=X$.
 Hence:
 \eq{ \Xnor = \Xlin + \sum_{n\in A(X)}B_n - [\Xlin,V] -\dotsb;}
 We set $V=\sum_{\mot{p}\in A(X)^*}V^{\mot{p}}B_{\mot{p}}$, then:
\eq{\Xnor = \Xlin + \sum_{n\in A(X)}B_n - \sum_{\mot{n}\in A(X)^*}\omega(\norme{\mot{n}}) V^{\mot{n}}
B_{\mot{n}}-\dotsb,}
after lemma~\ref{lem:xlin_1}. Hence the given expression of $\moule{V}$.

\end{proof}
\begin{theo}
\label{theo:simp}
 Let $V$ be the vector field $V = \sum_\bullet\moule{V}\comoule{B}$
where $\moule{V}$ is the mould given by:
\eq{ V^{\mot{n}} = \begin{cases}
                    \frac{1}{\omega(n)} &\text{ if $\longueur{n}=1$
                    and $\omega(n)\neq0$}\\
                    0 &\text{ otherwise.}
                   \end{cases}
}
  We call simplified form of $X$ the vector field $\Xsam =
\exp\bigl(V\bigr)X\exp\bigl(-V\bigr)$; it writes:
 \eq{ \Xsam = \Xlin + \sum_{\mot{n}\in A(X)^*}\Sam^{\mot{n}}
              B_{\mot{n}},}
and the mould $\moule{\Sam}$ has the following expression:
  \eq{
  \Sam^{\mot{n}} = \begin{cases}
                    0 &\text{ if $\mot{n}=\emptyset$,}\\
                    0 &\text{ if $\longueur{\mot{n}}=1$ and
                    $\omega(n)\neq0$,}\\
                    1 &\text{ if $\longueur{\mot{n}}=1$ and
                    $\omega(n)=0$;}
                   \end{cases}}
for the other words, the mould $\moule{\Sam}$ is given by the
following equation:
 \eq{ \moule{\Sam} = \moulexp{V}\times\moule{I}\times\moulexp{(-V)} -
 \nabla\moulexp{V}\times\moulexp{(-V)}.}
\end{theo}
\begin{proof}
The field $X=\Xlin +\sum_\bullet\moule{I}\comoule{B}$ has
non-resonant terms only of length $1$. If we want them to vanish, we
look for a simplified field $\Xsam = \Xlin +
\sum_\bullet\moule{\Sam}\comoule{B}$, where we set the mould
$\moule{\Sam}$ as follows:
 \eq{
  \Sam^{\mot{n}} = \begin{cases}
                    0 &\text{ if $\mot{n}=\emptyset$,}\\
                    0 &\text{ if $\longueur{\mot{n}}=1$ and
                    $\omega(n)\neq0$,}\\
                    &\text{we do not know yet for other words
                    $\mot{n}$.}
                   \end{cases}}
We keep the same equation as~\eqref{eq:conjugaison_3}:
 \eq{
  \moulexp{V}\times\moule{I} = \nabla \moulexp{V} +
  \moule{\Sam}\times\moulexp{V}.}
 By setting moreover that $V^{\mot{n}}=0$ if
 $\longueur{\mot{n}}\neq1$ it is possible to solve this equation,
 and find that:
  \eq{V^n = \begin{cases}
            0 &\text{ if $\omega(n)=0$,}\\
            \frac{1}{\omega(n)} &\text{ otherwise.}
            \end{cases}
            }
  Hence the simplified vector field $\Xsam$:
  \eqn{\label{eq:sam_1}\Xsam = \exp\Biggl( \sum_{\substack{n\in A(X)\\ \omega(n)\neq0}}
  \frac{1}{\omega(n)} B_n \Biggr) X \exp\Biggl( \sum_{\substack{n\in A(X)\\ \omega(n)\neq0}}
  -\frac{1}{\omega(n)}B_n \Biggr).}
\end{proof}
  The mould $\moule{\Sam}$ is then calculable, and we recall here
  its expression (see~\cite{calcul_moulien}):
    \begin{lemm}\label{lem:sam}
 The mould $\moule{\Sam}$ is given by:
 \begin{itemize}
 \item $\Sam^\emptyset=0$;
 \item if $\longueur{n}=1$ and $\omega(n)\neq0$, $\Sam^n=0$
   (kills the non-resonant terms);
 \item if $\longueur{n}=1$ and $\omega(s)=0$, $\Sam^n=1$;
 \item if $r=\longueur{\mot{n}}\geq2$ and
       $ \omega_1,\dotsc,\omega_r$ are different from $0$,
       \eq{\Sam^{\mot{n}} = \frac{1}{\omega_1\dotsm\omega_r}
              \sum_{k=1}^r
              \frac{(-1)^{r-k}(\omega_k(r-k)-
               \omega_{k+1}-\dotsb-\omega_r)}{(k-1)!(r-k+1)!}.}
 \item If only one $\omega_i$ vanishes,
       \eq{\Sam^{\mot{n}} = \frac{(-1)^{r-1}} {(i-1)!(r-i)!\omega_1\dotsm\omega_{i-1}
       \omega_{i+1}\dotsm\omega_r},}
 \item If more than one $\omega_i$ vanishes, then
       $\Sam^{\mot{n}}=0$.
 \end{itemize}
\end{lemm}
\subsection{Proof of lemma~\ref{lem:sam}}
 Let us denote by $\moule{V}$ the mould defined by:
 \eq{ V^{\mot{n}} = \begin{cases}
                        \frac{1}{\omega(n)} &\text{ if $\longueur{\mot{n}}=1$ and $\omega(n)\neq0$}\\
                        0 & \text{ otherwise}
                        \end{cases}.
                        }
  \subsubsection{A first term...}\label{subsub:C}
 We denote by $\moule{C}$ the mould $\moule{C}=\exp(\moule{V})\times\moule{I}$.
 We have $C^\emptyset=0$. For a word $n$ of length $1$ we have:
 \eq{C^n = \left(\exp(\moule{V})\right)^\emptyset I^n =1.}
 For a word $\mot{n}$ of length $r\geq2$,
  \eq{C^{\mot{n}} = \bigl(
      \exp(\moule{V})\bigr)^{n^1,\dotsc,n^{r-1}}I^{n^r}
      = \begin{cases}
        0  &\text{ if at least one of the $(\omega_i)_{1\leq i\leq r-1}$ vanishes,}\\
       \frac{1}{(r-1)!\omega_1\dotsm\omega_{r-1}} &\text{ otherwise.}
      \end{cases}
      }
  \subsubsection{...a second term...}
 We denote by $\moule{D}=\moule{C}\times\exp(-\moule{V})$, so
$D^\emptyset=0$, and for a word $n$ of length $1$, $D^n=C^n=1$.

  For a word $\mot{n}$ of length $r\geq2$ we have:
 \eq{D^{n^1,\dotsc,n^r} =
    C^{n^1}(\exp(-\moule{V}))^{n^2,\dotsc,n^r} +
    C^{n^1,n^2}(\exp(-\moule{V}))^{n^3,\dotsc,n^r} + \dotsb +
    C^{n^1,\dotsc,n^r}.}
 There are then several cases:
\begin{itemize}
   \item
  If at least one $\omega_i$ is zero, $1\leq i\leq r-1$,
       then every $C^{n^1,\dotsc,n^j}$, with $j\geq i+1$ vanishes
       (after the calculus of $C^{\mot{n}}$ in~\ref{subsub:C}); also every
       $(\exp(-\moule{V}))^{n^k,\dotsc,n^r}$ vanishes for $k\leq i$.
        Therefore $D^{\mot{n}}=C^{n^1,\dotsc,n^{i}}
       (\exp(-\moule{V}))^{n^{i+1},\dotsc,n^r}$.

        We have then for a word $\mot{n}$:
       \eq{D^{\mot{n}} = \begin{cases}
            0 &\text{ if another $(\omega_l)_{\substack{1\leq l\leq
            r\\l\neq i}}$ vanishes;}\\
            \frac{1}{(i-1)!\omega_1\dotsm\omega_{i-1}} \times
             \frac{(-1)^{r-i}}{(r-i)!\omega_{i+1}\dotsm\omega_r}
             &\text{ if no other $\omega_l, l\neq i$ is zero.}
           \end{cases}
           }
 \item if $\omega_r$ vanishes, $D^{\mot{n}}=C^{n^1,\dotsc,n^r}$
       therefore:
       \eq{D^{\mot{n}} =
           \begin{cases}
            0 &\text{ if one of the $(\omega_l)_{1\leq l\leq r-1}$
             vanishes;} \\
            \displaystyle\frac{1}{(r-1)!\omega_1\dotsm\omega_{r-1}}
              &\text{ otherwise.}
           \end{cases}
           }
 \item  if no other $\omega_i$ vanishes then
       \begin{multline*}
        D^{n^1,\dotsc,n^r} = \frac{(-1)^{r-1}}{(r-1)!\omega_2\dotsm\omega_r} +
        \frac{1}{\omega_1}\times\frac{(-1)^{r-2}}{(r-2)!\omega_3\dotsm\omega_r}+
        \dotsb\\
        \qquad+\frac{1}{(r-2)!\omega_1\dotsm\omega_{r-2}}\times\frac{-1}{\omega_r}
        +\frac{1}{(r-1)!\omega_1\dotsm\omega_{r-1}},
       \end{multline*}
       that is:
        \eq{D^{n^1,\dotsc,n^r}=\frac{1}{\omega_1\dotsm\omega_r}
            \sum_{k=1}^r\frac{(-1)^{r-k}\omega_k}{(k-1)!(r-k)!}.}
\end{itemize}
 \subsubsection{...a third term...}
  We calculate now the following:
 $\moule{E}=\nabla\exp(-\moule{V})$. Thanks to the previous computations on
 the mould exponential, and by definition of $\nabla$, we have $E^\emptyset=0$;
  moreover, for a word $\mot{n}$ of length $r\geq1$:
  \eq{E^{n^1,\dotsc,n^r} = (\omega_1+\dotsb+\omega_r)(\exp(-\moule{V}))^{\mot{n}},}
 then
  \eq{E^{\mot{n}} =
      \begin{cases}
       0 &\text{ if one at least of the $(\omega_i)_{1\leq i\leq r}$
         vanishes;}\\
       \frac{(\omega_1+\dotsb+\omega_r)(-1)^r}{r!\omega_1\dotsm\omega_r} &\text{ otherwise.}
      \end{cases}
      }
 \subsubsection{... the last term}
We finally calculate the mould product
  $\moule{F}=\exp(\moule{V})\times\moule{E}$~; first,
$F^\emptyset=0$~; for a word of length $1$, $F^n=E^n$ therefore
$F^n=0$ if $\omega(n)=0$ and $F^n=-1$ if $\omega(n)\neq0$~; at last,
for a word $\mot{n}$ of length $r\geq1$ we have:
 \eq{F^{n^1,\dotsc,n^r} = (\exp(\moule{V}))^\emptyset E^{\mot{n}} +
      (\exp(\moule{V}))^{n^1}E^{n^2,\dotsc,n^r} + \dotsb +
      (\exp(\moule{V}))^{n^1,\dotsc,n^{r-1}}E^{n^r}.}
 Thus, after the calculus of $\exp(\moule{V})$ and $\moule{E}$,
$F^{\mot{n}}$ vanishes if at least one of the $\omega_i$ is zero.
 If no $\omega_i$ is zero, then:
 \begin{align*}
 F^{n^1,\dotsc,n^r} &=
 \frac{(\omega_1+\dotsb+\omega_r)(-1)^r} {r!\omega_1\dotsm\omega_r}
 + \frac{1}{\omega_1} \times \frac{(\omega_2+\dotsb+\omega_r)(-1)^{r-1}}
 {(r-1)!\omega_2\dotsm\omega_r} +\dotsb\\
 &\quad +\frac{1}{(r-1)!\omega_1\dotsm\omega_{r-1}} \times
\frac{(-1)\omega_r}{\omega_r}.
 \end{align*}
Finally:
 \eq{F^{\mot{n}} =
     \begin{cases}
      0 &\text{ if one of the $\omega_i$ is zero,}\\
      \frac{1}{\omega_1\dotsm\omega_r}
       \sum_{k=1}^{r}\frac{(-1)^{r-k+1}(\omega_{k}+\dotsb+\omega_r)}
       {(r-k+1)!(k-1)!} &\text{ otherwise.}
     \end{cases}
     }
 \subsubsection{Expression of the sought mould}
 As $\moule{\Sam}=\moule{F}+\moule{D}$, we have
 \eq{\Sam^\emptyset = F^\emptyset+D^\emptyset = 0.}
For a word of length $1$, we get:
  \eq{\Sam^n =
      \begin{cases}
       1 &\text{ if $\omega(n)=0$,}\\
       0 &\text{ if $\omega(n)\neq0.$}
      \end{cases}
      }
  The non-resonant terms of the field $X$ thus vanish in the field $\Xsam$.
 Moreover, for a word $\mot{n}$ of length $r\geq2$, we get:
 \begin{itemize}
 \item if there exists $i$ and $j$, two different integers from
       $\{1,\dotsc,r\}$ such as $\omega_i=\omega_j=0$ then
       $D^{\mot{n}}=F^{\mot{n}}=0$ therefore
       \eq{\Sam^{\mot{n}}=0.}
 \item if one $\omega_i$ exactly is zero, $F^{\mot{n}}=0$ and
       \eq{\Sam^{\mot{n}} = \frac{1}{(i-1)!\omega_1\dotsm\omega_{i-1}} \cdot
           \frac{(-1)^{r-1}}{(r-i)!\omega_{i+1}\dotsm\omega_r}.}
 \item if every $(\omega_i)_{1\leq i\leq r}$ is  non zero,
       then \eq{D^{\mot{n}}=\frac{1}{\omega_1\dotsm\omega_r}
                 \sum_{k=1}^r\frac{(-1)^{r-k}\omega_k}{(k-1)!(r-k)!},}
             \eq{F^{\mot{n}} = \frac{1}{\omega_1\dotsm\omega_r}
                 \sum_{k=1}^{r}\frac{(-1)^{r-k+1}(\omega_{k}+\dotsb+\omega_r)}
                 {(r-k+1)!(k-1)!},}
       and \eq{\Sam^{\mot{n}} = \frac{1}{\omega_1\dotsm\omega_r}
              \sum_{k=1}^r
              \frac{(-1)^{r-k}(\omega_k(r-k)-\omega_{k+1}-\dotsb\omega_r)}{(k-1)!(r-k+1)!}.
              }
 \end{itemize}
The proof is now complete!

\input{rajout}
 The field $\Xsam$ has now no more non-resonant terms of length $1$, but the
 transformation $X \mapsto \Xsam$, that we denote $\simp$, introduces
 non-resonant terms of length greater than $1$, as $\Sam^{\mot{n}}$ is
 not always $0$ when $\longueur{\mot{n}}\geq2$.
\subsection{The first step}
  Now that we have killed every non-resonant terms of length $1$,
  but introduced some more of length greater than $2$, we have
  to iterate the transform... however, if we want the iteration procedure to be writable
  in terms of moulds composition, we
  \emph{must} find a way to rewrite $\Xsam = \Xlin + \sum_{\mot{n}\in A(X)^*}
  \Sam^{\mot{n}} B_{\mot{n}}$ as $\Xsam = \Xlin + \sum_{m\in\A}D_m$
  where $\A$ is a new alphabet and $(D_m)_{m\in\A}$ a new collection
  of homogeneous differential operators. It is here natural at
  this time, since for every word $\mot{n}$ in $A(X)^*$,
  $B_{\mot{n}}$ is a homogeneous differential operator of degree
  $\norme{\mot{n}}$ (nevertheless, do not forget that the most
  natural choice may not always be the best, as we will see in
  section~\ref{section:hamilton}).
  Hence the new alphabet is:
  \eq{\A=A\bigl(\Xsam^{(1)}\bigr) = \bigl\{ \norme{\mot{n}},
        \mot{n} \in A(X)^* \bigr\} \subset
  \Z^\nu.}
Let us write then $\Xsam^{(1)}$ instead of $\Xsam$, for it is the
first of the iteration:
\eq{ \Xsam^{(1)} = \Xlin + \sum_{m\in\A} D_m,
    \text{ with $D_m = \sum_{\substack{\mot{n}\in A(X)^*\\
                             \norme{\mot{n}}=m}} \Sam^{\mot{n}}
                             B_{\mot{n}}$,}}
and do the transformation $\simp$ again. We get a $\Xsam^{(2)}$ and:
 \eq{ \Xsam^{(2)} = \Xlin + \sum_{\mot{m}\in\A^*} \Sam^{\mot{m}}
      D_{\mot{m}.}}
Writing this in the old alphabet, we have, by definition of the
composition of two moulds:
 \eq{ \Xsam^{(2)} = \Xlin + \sum_{\mot{n}\in A(X)^*} \bigl(
 \moule{\Sam}\circ\moule{\Sam} \bigr)^{\mot{n}} B_{\mot{n}}.}

%% file: rajout.tex
In order to put in evidence the universal feature of the moulds $\Sam$ we
obtain under the simplification procedure, we introduce the
following one parameter family of complex valued functions:

\begin{defi}
The \emph{Poincaré family} is denoted by $\cal{P}=\left ( P_q \right)_{q\in \N^*}$ where $P_q :\C^q \rightarrow \C^q$, and defined by $P_1 (z)=1$ if $z=0$ and $0$ otherwise, and for $q\geq 2$
\eq{
 P_q (z_1 ,\dots ,z_q )=
  \begin{cases}
 \frac{1}{z_1\dotsm z_q} \sum_{k=1}^q  (-1)^{q-k}  \frac{z_k (q-k) -z_{k+1} - \dotsm - z_q }{(k-1)! (q-k+1)!},  &\text{ if $z\in (\C^* )^q$ ,}  \\
 \frac{(-1)^{q-1}}{(i-1)! (q-i)!}  \frac{1}{ z_1 \dots z_{i-1} z_{i+1} \dots z_q}, &\text{ if $z\in S_{q,i}$},\\
 0 &\text{ otherwise,}
\end{cases}
}
with $S_{q,i}=(\C^* )^{i-1} \times \{ 0\} \times (\C^* )^{q-i}$.

\end{defi}

Lemma~\ref{lem:sam} of this section can then be formulated as follows:

\begin{lemm}
The simplification moulds $\moule{\Sam}$ is given by $\Sam^{\emptyset} =1$ and $\Sam^{\mot{n}} =P_{\longueur{\mot{n}}} (\mot{n} .\lambda )$ if
$\longueur{\mot{n}}\geq 1$.
\end{lemm}

%% file: trimmed.tex
%
\label{section:trimmed}
\subsection{The simplification procedure}
 We can construct a prenormal form by iterating the simplification
  procedure seen in the
 previous section.
  \begin{defi}[Trimmed form up to order $r$]
Let $r$ be in $\N$. The trimmed form up to order $r$ is defined as
$\Xsam^{(r)}$, obtained from $X$ after $r$ successive
simplifications:
  \eq{ X = \Xsam^{(0)} \xrightarrow{\simp_1} \Xsam^{(1)}
  \xrightarrow{\simp_2}\dotsb \xrightarrow{\simp_r}
  \Xsam^{(r)},}
where $\simp_i$ is the automorphism of simplification defined by:
 \eq{ \simp_i = \exp\bigl(V_i\bigr),}
with $V_i$ the vector field associated to the mould $\moule{V}$ on
the alphabet $A\bigl(\Xsam^{(i)}\bigr)$ defined recursively by:
 \eq{A\bigl(\Xsam^{(i)}\bigr) = \bigl\{ \norme{\mot{n}}, \mot{n} \in
                                    A\bigl(\Xsam^{(i-1)}\bigr)^* \bigr\}.}
\end{defi}
 Theorem~\ref{theo:simp} leads us then to the following result:
\begin{theo}
 For all $r$ in $\N$, the trimmed form up to order $r$ $\Xsam^{(r)}$
 has a mould expansion, \ie there exists a mould denoted by
 $\moule{\Sam}_r$ in $\mathcal{M}_{\C}\bigl(A(X)\bigr)$ such that:
  \eq{ \Xsam^{(r)} = \Xlin + \sum_{\bullet} \moule{\Sam}_r\comoule{B}.}
Moreover the mould $\moule{\Sam}_r$ can be defined with the help of
$\moule{\Sam}$:
 For all $r$ in $\N$, we have:
 \eq{ \moule{\Sam}_r = \underbrace{\moule{\Sam} \circ\dotsb\circ
 \moule{\Sam}}_{r \text{ times}}.}
\end{theo}
 From now on, we denote either $\moule{\Sam}_r$ or $\bigl(\moule{\Sam}\bigr)^{\circ r}$
  for the composition of $r$ copies of $\moule{\Sam}$.
Let us investigate now what is happening: we already saw that after
the first step, non-resonant terms of length $1$ --actually, there
are no others-- vanished. We have then:
\eq{ \Xsam^{(1)} = \Xlin + \sum_{\substack{n\in A(X)\\ \omega(n)=0}}
    B_n + \sum_{m\in A(X)} \underbrace{\sum_{\substack{\mot{n}\in A(X)^*\\
    \norme{\mot{n}}=m}} \Sam^{\mot{n}} B_{\mot{n}}}_{D_m}.}
Now we see from such a writing, that the simplification introduces
new terms, which may not be (and usually are not) resonant.

  The main property of $\moule{\Sam}_r$ is that it is "stationary"
 in this sense:
\begin{theo}\label{theo:stationaire}
 Let $r$ be in $\N^*$.
For any word $\mot{n}$ in $A^*$ of length at most $r$, we have:
\eq{\Sam_q^{\mot{n}} = \Sam_r^{\mot{n}},\,\forall q\geq r.}
\end{theo}
This theorem is deduced from the following lemma:
\begin{lemm}
Let $r$ be in $\N^*$. For any word $\mot{n}$ in $A^*$ of length at
most $r$, we have:
\eq{\Sam_{r+1}^{\mot{n}} = \Sam_r^{\mot{n}}.}
\end{lemm}

\begin{proof}[Proof of the lemma]
The proof is done by induction on the length $r$.

 For $r=1$: let $n$ be in $A^*$, we have:
\eq{ \bigl(\moule{\Sam}\circ\moule{\Sam}\bigr)^n =
\bigl(\Sam^n\bigr)^2 = \Sam^n, \text{ after lemma~\ref{lem:sam}}.}
 We suppose now that $r\geq2$ is fixed, and that for any $p\leq
 r-1$, and any word $\mot{n}$ of length at most $p$, $\Sam_{p+1}^{\mot{n}}
  = \Sam_p^{\mot{n}}$. Let $\mot{n}=n^1\dotsb n^r$ be a word of length $r$. We
have three cases:
\begin{enumerate}
 \item
 $\mot{n}$ is non-resonant, \ie $\omega(\norme{\mot{n}})\neq0$. In
 this case we write:
\eq{
 \begin{aligned}
 \Sam_{r+1}^{\mot{n}} &= \bigl(\moule{\Sam}\circ
 \moule{\Sam}_r\bigr)^{\mot{n}} \\
  &= \Sam^{\norme{\mot{n}}} \Sam_r^{\mot{n}} + \sum_{\substack{
  2\leq l\leq r\\ \mot{n}_1\dotsb\mot{n}_l=\mot{n}}}
  \Sam^{\norme{\mot{n}_1},\dotsb,\norme{\mot{n}_l}}
  \Sam_r^{\mot{n}_1}\dotsm \Sam_r^{\mot{n}_l}.
  \end{aligned}
  }
 As $\omega(\norme{\mot{n}})\neq0$, by lemma~\ref{lem:sam},
 $\Sam^{\norme{\mot{n}}}=0$, hence
  $ \Sam^{\norme{\mot{n}}} \Sam_r^{\mot{n}}
   = \Sam^{\norme{\mot{n}}} \Sam_{r-1}^{\mot{n}}$.
 Moreover, in the second term, as $l\geq2$,
 $\longueur{\mot{n}_k}\leq r-1$ for any $k$ in $\{ 1,\dotsc,l\}$,
 hence by the induction assumption,
 $\Sam_r^{\mot{n}_k}=\Sam_{r-1}^{\mot{n}_k}$.
 Finally,
 \eq{ \begin{aligned}
     \Sam_{r+1}^{\mot{n}} &=
     \Sam^{\norme{\mot{n}}}\Sam_{r-1}^{\mot{n}} + \sum_{\substack{
  2\leq l\leq r\\ \mot{n}_1\dotsb\mot{n}_l=\mot{n}}}
  \Sam^{\norme{\mot{n}_1},\dotsb,\norme{\mot{n}_l}}
  \Sam_{r-1}^{\mot{n}_1}\dotsm \Sam_{r-1}^{\mot{n}_l} \\
   &= \bigl(\moule{\Sam}\circ
 \moule{\Sam}_{r-1}\bigr)^{\mot{n}} \\
   &= \Sam_r^{\mot{n}}.
   \end{aligned}
   }
 \item
 $\mot{n}$ is resonant, \ie $\omega(\norme{\mot{n}})=0$, and
 $\omega(n^i)=0$ for all $i$ in $\{1,\dotsc,r\}$. In this case, we
 write again:
 \eq{
 \begin{aligned}
  \Sam_{r+1}^{\mot{n}} &= \bigl(\moule{\Sam}\circ
 \moule{\Sam}_r\bigr)^{\mot{n}} \\
  &= \Sam^{\norme{\mot{n}}} \Sam_r^{\mot{n}} + \sum_{\substack{
  2\leq l\leq r\\ \mot{n}_1\dotsb\mot{n}_l=\mot{n}}}
  \Sam^{\norme{\mot{n}_1},\dotsb,\norme{\mot{n}_l}}
  \Sam_r^{\mot{n}_1}\dotsm \Sam_r^{\mot{n}_l}.
  \end{aligned}
}
 On the one hand, after lemma~\ref{lem:sam}, $\Sam^{\norme{\mot{n}}}=1$, for
  $\omega(\norme{\mot{n}})=0$; on the other hand, for $l$ in
  $\{2,\dotsc,r\}$, $\omega(\norme{\mot{n}_k})=0$ for all $k$ in $\{
  1,\dotsc,l\}$, hence, still after lemma~\ref{lem:sam},
  $\Sam^{\norme{\mot{n}_1},\dotsb,\norme{\mot{n}_l}} = 0$.
  Finally, $\Sam_{r+1}^{\mot{n}}= \Sam_r^{\mot{n}}$.
 \item
 $\mot{n}$ is resonant, \ie $\omega(\norme{\mot{n}})=0$ and there is
 at least one (therefore two) $n^i$ in $\mot{n}$ such that
 $\omega(n^i)\neq0$. After lemma~\ref{lem:sam}, $\Sam^{n^i}=0$.
  In this last case, we write:
 \eq{ \begin{aligned}
 \Sam_{r+1}^{\mot{n}} &= \bigl( \moule{\Sam}_r \circ
 \moule{\Sam}\bigr)^{\mot{n}} \\
  &= \sum_{\substack{1\leq l\leq r-1 \\
  \mot{n}_1\dotsb\mot{n}_l=\mot{n}}}
  \Sam_r^{\norme{\mot{n}_1}\dotsb\norme{\mot{n}_l}} \Sam^{\mot{n}_1}
  \dotsm \Sam^{\mot{n}_l} + \Sam_r^{\mot{n}}
  \Sam^{n^1}\dotsm\Sam^{n^r}.
  \end{aligned}
  }
 By the induction assumption, for every $l$ in $\{1,\dotsc,r-1\}$,
  \eq{ \Sam_r^{\norme{\mot{n}_1}\dotsb\norme{\mot{n}_l}} =
   \Sam_{r-1}^{\norme{\mot{n}_1}\dotsb\norme{\mot{n}_l}};}
moreover:
\eq{
\sum_{\substack{1\leq l\leq r-1 \\
  \mot{n}_1\dotsb\mot{n}_l=\mot{n}}}
  \Sam_{r-1}^{\norme{\mot{n}_1}\dotsb\norme{\mot{n}_l}}
   \Sam^{\mot{n}_1}
  \dotsm \Sam^{\mot{n}_l}
   = \bigl( \moule{\Sam}_{r-1}\circ\moule{\Sam} \bigr)^{\mot{n}} -
   \Sam_{r-1}^{\mot{n}} \Sam^{n^1}\dotsm\Sam^{n^r}.}
Hence
 $ \Sam_{r+1}^{\mot{n}} = \Sam_r^{\mot{n}} -
  \Sam_{r-1}^{\mot{n}} \Sam^{n^1}\dotsm\Sam^{n^r} + \Sam_r^{\mot{n}}
  \Sam^{n^1}\dotsm\Sam^{n^r}$.
Now, after lemma~\ref{lem:sam}, the product
$\Sam^{n^1}\dotsm\Sam^{n^r}$ is zero.
 Finally, $\Sam_{r+1}^{\mot{n}} = \Sam_r^{\mot{n}}$.
\end{enumerate}
\end{proof}
\subsection{The Poincaré-Dulac theorem}
 We can define now the mould $\moule{\Tram}$ as follows:
\begin{defi}
\label{defi:tram} The mould $\moule{\Tram}$ is defined by
$\Tram^\emptyset=0$, and for a word $\mot{n}$ in $A^*$ of length
$r\geq1$,
 \eq{ \Tram^{\mot{n}} = \Sam_r^{\mot{n}} = \lim_{p\to+\infty}
 \bigl(\bigl( \moule{\Sam} \bigr)^{\circ p} \bigr)^{\mot{n}}
 }
The limit exists after theorem~\ref{theo:stationaire}.
\end{defi}
 We define then the trimmed form.
 \begin{defi}
  The trimmed form of $X$ is the limit of the simplification
  procedure. It is given by:
  \eq{ \Xtram = \Xlin + \sum_{\bullet} \moule{\Tram}\comoule{B}.}
 \end{defi}
 Now this result shows that the trimmed form is what we are looking for:
\begin{theo}
\label{theo:tram}
 The trimmed form is a continuous prenormal form.
 \end{theo}
\begin{proof}
 Remember that $\Tram^\emptyset=0$. By definition of
 $\moule{\Tram}$, for a word of length $1$, $\Tram^{n}=\Sam^{n}$,
 hence $\Tram^n=0$ if $\omega(n)\neq0$.
 Now, for a word $\mot{n}$ of length greater than $2$, by
 definition~\ref{defi:tram} of $\moule{\Tram}$ we have:
\begin{align}
 \moule{\Tram} &= \moule{\Tram}\circ\moule{\Sam}\label{eq:tram_1}\\
             &= \moule{\Sam}\circ\moule{\Tram}.\label{eq:tram_2}
\end{align}
We can then verify that $\Tram^{\mot{n}}=0$ if
$\omega(\norme{\mot{n}})\neq0$ by induction on the length $r\geq2$
of $\mot{n}$. If $\mot{n}=n^1n^2$ is in $A^*$, by definition of the
composition of two moulds and after~\eqref{eq:tram_1}:
 \eq{\begin{aligned}
     \Tram^{\mot{n}} &= \Tram^{\emptyset}\Sam^{\mot{n}} + \Tram^{n^1}\Sam^{n^2}
                      + \Tram^{\mot{n}}\Sam^{\emptyset} \\
              &=  \Tram^{n^1}\Sam^{n^2}\\
              &= \Sam^{n^1}\Sam^{n^2}.
     \end{aligned}
     }
Now, after lemma~\ref{lem:sam}, $\Tram^{\mot{n}}\neq0$ if and only
if $\omega(n^1)=\omega(n^2)=0$, and this is impossible since
$\omega(n^1)+\omega(n^2)\neq0$. Hence $\Tram^{\mot{n}}=0$ if
$\omega(\norme{\mot{n}})\neq0$.

 We fix $r\geq3$ and suppose that $\Tram^{\mot{n}}=0$ if
$\omega(\norme{\mot{n}})\neq0$, for any word $\mot{n}$ of length
less than $r-1$.
 Then if $\mot{n}$ is a word of length $r$ such that
  $\omega(\norme{\mot{n}})\neq0$ we have after equation~\eqref{eq:tram_2}:
 \eq{
     \Tram^{\mot{n}}  = \Sam^{\norme{\mot{n}}}\Tram^{\mot{n}} +
     \sum_{\substack{2\leq l\leq r\\
     \mot{n}_1\dotsb\mot{n}_l=\mot{n}}}
     \Sam^{\norme{\mot{n}_1}\dotsb\norme{\mot{n}_l}}
     \Tram^{\mot{n}_1}\dotsm\Tram^{\mot{n}_l}.}
 The term $\Sam^{\norme{\mot{n}}}\Tram^{\mot{n}}$ is $0$, for
 $\omega(\norme{\mot{n}})\neq0$. Now for each partition of
 $\mot{n}$ in $l$ words $\mot{n}_1,\dotsc,\mot{n}_l$, where $l\geq2$
 there is at least one $k$ in $\{1,\dotsc,l\}$ such that
 $\omega(\norme{\mot{n}_k})\neq0$ (for $\omega(\norme{\mot{n}}) =
 \sum_{j=1}^l \omega(\norme{\mot{n}_j})$). Hence by induction,
 $\Tram^{\mot{n}_k}=0$. Therefore, $\Tram^{\mot{n}}=0$.
 \end{proof}

%% file: hamiltonian.tex
\label{section:hamilton}
 We discuss here the application of the preceding sections to Hamiltonian
 operators. $H$ is a Hamiltonian function, in cartesian coordinates:
 \eq{ H(x,y)= \sum_{i=1}^\nu \lambda_i x_i y_i +
 \sum_{(n,m)\in A(H)} a_{nm} x^n y^m,}
 where $A(H)$ stands for the set of higher degrees and
 $\mot{\lambda}=(\lambda_1,\dotsc,\lambda_\nu)$ an element of
 $\C^\nu$. We will denote $A$ when there is no ambiguity.
 The Hamiltonian vector field then writes:
 \eq{
  X_H = \Xlin + \sum_{(n,m)\in A}
  D_{nm},}
  where
  \eq{ \Xlin = -\sum_{i=1}^{\nu}\lam_i x_i \dep_{x_i}+\sum_{i=1}^{\nu}\lam_i
  y_i\dep_{y_i}}
and
   \eq{D_{nm} =a_{nm}\sum_{i=1}^{\nu} x^{\hat{n}_i}y^{\hat{m}_i}(n_i
  y_i\dep_{y_i} - m_i x_i\dep_{x_i}),}
  and we denote $\hat{n}_i = (n_1,\dots,n_i-1,\dots,n_\nu)$ (same for
  $m_i$).
   Remark that $D_{nm}$ is \emph{not} an homogeneous operator; however $D_{nmi}$
    defined as follows is a homogeneous operator of degree $(\hat{n}_i,\hat{m}_i)$:
    \eqn{D_{nmi} = a_{nm}x^{\hat{n}_i}y^{\hat{m}_i}(n_i
  y_i\dep_{y_i} - m_i x_i\dep_{x_i}),}
  and $D_{nm} = \sum_{i=1}^\nu D_{nmi}$.
 The preceding "usual" decomposition  in homogeneous operators of section~\ref{subsec:resonant}
  does not lead to
 Hamiltonian operators. Nevertheless lemma~\ref{lem:proj_hamilton}
gives a way to obtain Hamiltonian operators... when starting also
from Hamiltonian ones.

 The interesting thing about $D_{nm}$ is that it is a Hamiltonian
 operator, \ie it defines a Hamiltonian vector field. We will
 frequently denote $s=(n,m)$ a letter of $A$, and $D_s$ for
 $D_{nm}$. As previously, if $\mot{s}=s^1\dotsb s^r$ is a word in
 $A^*$, $D_{\mot{s}}$ will be the composition $D_{s^1}\dotsm
 D_{s^r}$.
 It is well-known that a prenormal form of a Hamiltonian vector field
 is also a Hamiltonian vector field, and that the transformation which
 brings the former into the latter is symplectic.
 
 However, it is also important to keep in mind that we want \emph{successive}
 canonical transformations to preserve the Hamiltonian character, because, for
 example, if we want to implement that prenormal form, a computer cannot do an infinite number
 of iterations.
 
 Nevertheless, if we decompose $X_H$ in homogeneous differential
 operators, as done before, it is very difficult to know if we get Hamiltonian transformation!

 The following lemma gives a first result on "Hamiltonian-preserving" moulds and justifies the
 use of the $D_{nm}$ operators instead of usual homogeneous operators. We will
 need this result in the next subsection.
\begin{lemm}\label{lem:proj_hamilton}
 Let $\moule{M}$ be an alternal mould on an alphabet $A^*$. Let
 $\comoule{S}$ be a collection of differential operators, such that
 $S_u$ is a Hamiltonian vector field for every $u$ in $A$. Then the
 sum $\sum_{\mot{u}\in A^*} M^{\mot{u}} S_{\mot{u}}$
  defines a Hamiltonian vector field.
\end{lemm}
\begin{proof}
 The key is that if $S_u=X_{H_u}$ and $S_t=X_{H_t}$ are Hamiltonian
 vector fields, then:
  \eq{ [ S_u,S_t ] = \{ X_{H_u},X_{H_t} \} = X_{\{H_u,H_t\}},}
  where $\{,\}$ is the usual Poisson bracket; hence $[S_u,S_t]$ is
  still Hamiltonian. Thus, by an induction on $r$, for any word $\mot{u}=u^1\dotsb u^r$ of
$A^*$ of length $r$, $S_{[\mot{u}]}$ is a Hamiltonian vector field.
    We now use the projection lemma~\ref{lem:projection}: $\moule{M}$ being
alternal, we have, if $\mot{u}$ is a word of length $r$ and
$\sigma(\mot{u})$ the set of words deduced from $\mot{u}$ by a
permutation:
 \eq{ \sum_{\mot{u}\in \sigma(\mot{u})} M^{\mot{u}} S_{\mot{u}} =\frac{1}{r}
 \sum_{\mot{u}\in \sigma(\mot{u})} M^{\mot{u}} S_{[\mot{u}]}.}
  Let us denote $\sim$ the equivalence relation on $\alphabet{A}$
    defined by:
    \eq{ \mot{u}\sim\mot{t} \iff \text{there exists one permutation
    $\tau$ such as $\tau(\mot{u}) = \mot{t}$.}}
    We have $\mot{u}\sim\mot{t} \iff \longueur{\mot{u}} = \longueur{\mot{t}}$
    therefore $\alphabet{A}\bigl\backslash \{\emptyset\}\bigr. = \coprod_{r=1}^{+\infty}
    \alphabet{A_r}$ and $\alphabet{A}\Bigl/_{\!\displaystyle\sim}\Bigr. =
    \coprod_{r=1}^{+\infty} \alphabet{A_r}\Bigl/_{\!\displaystyle\sim}\Bigr.$; moreover
    $M^\emptyset=0$, for $\moule{M}$ is alternal; hence the
    following equalities:
    \eq{
     \begin{aligned}
      \sum_{\mot{u}\in\alphabet{A}} M^{\mot{u}}S_{\mot{u}}
               &= \sum_{r\geq1}\sum_{\mot{u}\in\alphabet{A_r}}
               M^{\mot{u}} S_{\mot{u}} \\
               &= \sum_{r\geq1} \sum_{\mot{u}\in\alphabet{A_r}\bigl/_{\!\scriptstyle\sim}\bigr.}
                  \sum_{\mot{u}\in\sigma(\mot{u})} M^{\mot{u}}
                  S_{\mot{u}} \\
               &= \sum_{r\geq1}\sum_{\mot{u}\in\alphabet{A_r}\bigl/_{\!\scriptstyle\sim}\bigr.}
                  \frac{1}{r} \sum_{\mot{u}\in\sigma(\mot{u})}
                  M^{\mot{u}} S_{[\mot{u}]} \\
               &= \sum_{r\geq1} \frac{1}{r}
               \sum_{\mot{u}\in\alphabet{A_r}} M^{\mot{u}}
               S_{[\mot{u}]}.
      \end{aligned}
     }
 Hence the result.
\end{proof}
 From now on, for $s=(n,m)$ in $A(H)$, we denote
 $\omega(s)=\omega(n,m) = \sum_{j=1}^\nu \lambda_j(m_j-n_j)$.
 As previously, a word $\mot{s}$ is \emph{resonant} if $\omega(\norme{\mot{s}})=0$.
  We have an analogous result as lemma~\ref{lem:xlin_1}:
 \begin{lemm}\label{lem:xlin_2}
For $\mot{s}=s^1\dotsb s^r$ a word in $A^*$ of length $r$, we have:
 \eq{ \Xlin D_{\mot{s}} = D_{\mot{s}} \Xlin + \norme{\omega(\mot{s})} D_{\mot{s}}.}
 \end{lemm}
\subsection{The limit of the simplification procedure}
 We proceed as before, by associating to $X_H$ a simplified vector
 field $\Xsam$ in the following way:
 \eq{\Xsam = \exp\Biggl( \sum_{\substack{s\in A(H)\\ \omega(s)\neq0}}
  \frac{1}{\omega(s)} D_s \Biggr) X_H \exp\Biggl( \sum_{\substack{s\in A(H)\\ \omega(s)\neq0}}
  -\frac{1}{\omega(s)}D_s \Biggr).}
 The important thing is that
\eq{\Xsam = \Xlin + \sum_{\mot{s}\in
 A(H)^*} \Sam^{\mot{s}} D_{\mot{s}},
}
 with $\moule{\Sam}$ \emph{exactly} the same mould as defined in
 lemma~\ref{lem:sam}. The only things that change are the alphabet
 and the operators. But the fact that $\moule{\Sam}$ is alternal
 is still true of course: $\Xsam$ is then, by
 lemma~\ref{lem:proj_hamilton}, a Hamiltonian vector field, and by
 definition of $\moule{\Sam}$ a trimmed form of $X_H$ up to order
 $1$.
 
 We then want to rewrite $\Xsam$ as a sum $\Xlin + \sum_{a\in\A}\Delta_a$
 where $\A$ would be a new alphabet, and there would be a simple law $\star$ such
 that $\norme{\mot{n}}_\star=a$ for a word $\mot{n}$ of $A(X)^*$. The second step
 would then be given by the composition $\moule{\Sam}\circ\moule{\Sam}$. Unfortunately,
 we have not been able, so far, to find such a new alphabet to make the iteration easy
 to formulate. So, we changed --again!-- the decomposition of the initial vector field $X_H$,
 therefore the alphabet, so that a mould iteration can be done.
\subsection{Canonical simplification}
 We define a new alphabet $\Omega(H)$, or $\Omega$ when there is no
 ambiguity, by:
 \eq{ \Omega(H)=\bigl\{ \omega(s) \text{ with } s \in A(H) \}.}
Remark that $\Omega$ is thus part of $\C$ and not anymore of
$\Z^\nu$.
     We have then
     \eq{ H(x,y) = \sum_{i=1}^\nu \lambda_ix_iy_i +
     \sum_{\omega\in\Omega} \sum_{\substack{(n,m)\in A(H)\\
 \omega(n,m)=\omega}} a_{nm}x^ny^m .}
\begin{defi}
 For $\omega\in\Omega$, $\D_\omega$ is the Hamiltonian
     vector field induced by the sum of monomials $H_\omega(x,y)=\sum_{\substack{(n,m)\in A(H)\\
     \omega(n,m) = \omega}} a_{nm} x^ny^m$; we call this latter sum
     the $\Omega$-homogeneous component of degree $\omega$ of $H$.
      We write
     $\D_\omega=X_{H_\omega}$.
\end{defi}
  The $\D_\omega$ are still Hamiltonian operators, as sum of such
  operators: in fact $\D_\omega = \sum_{\substack{(n,m)\in A(H)\\ \omega(n,m)=\omega}}
 D_{nm}$.
 We have thus $X_H=\Xlin + \sum_{\omega\in\Omega} \D_\omega$.
  This gives us the action of $\Xlin$ on the
 $\D_\omega$ (analogous to lemma~\ref{lem:xlin_2}):
\begin{lemm}
  \eq{ \Xlin \D_{\mot{\omega}} = \D_{\mot{\omega}}\Xlin + \norme{\mot{\omega}} \D_{\mot{\omega}}
    \text{ for any word $\mot{\omega}\in\alphabet{\Omega}$.}}
\end{lemm}
 The simplified field $\Xsam$ is obtained exactly the same way as
above:
  \eq{\Xsam = \Biggl( \sum_{\substack{\omega\in \Omega(H)\\
\omega\neq0}}
  \frac{1}{\omega} \D_\omega \Biggr) X_H \Biggl( \sum_{\substack{\omega\in \Omega(H)\\
    \omega\neq0}} -\frac{1}{\omega}\D_\omega \Biggr),}
and still:
 \eq{ \Xsam = \Xlin +
 \sum_{\mot{\omega}\in\Omega^*}\Sam^{\mot{\omega}}\D_{\mot{\omega}}.}
The mould $\moule{\Sam}$ is again defined as in lemma~\ref{lem:sam},
but on the alphabet $\Omega$, so for $\mot{\omega}=\omega(s^1)\dotsb
\omega(s^r)$ in $\Omega^*$, we set $\Sam^{\mot{\omega}} =
\Sam^{\mot{s}}$.
 $\moule{\Sam}$ is alternal, then $\Xsam$ is a Hamiltonian vector
 field.

 Now we want to iterate this process, as we did at the beginning of
section~\ref{section:trimmed}. In order to iterate, we must rewrite
$\Xsam$ as:
 \eq{ \Xsam = \Xlin + \sum_{\tilde{\omega}\in\tilde{\Omega}}
 \D^{(1)}_{\tilde{\omega}},}
where $\tilde{\Omega}$ is the new alphabet, and
 $\D^{(1)}_{\tilde{\omega}}$ expresses with the $\D_\omega$ and is
 still Hamiltonian. By definition, $\D_{\tilde{\omega}}^{(1)}$
  is the Hamiltonian field
coming from the sum $\sum_{\omega(n,m)=\tilde{\omega}}
a^{(1)}_{nm}x^ny^m$ in the new Hamiltonian. So, for a
$\tilde{\omega}$ fixed, we must find the $\mot{\omega}$ such that
$\D_{\mot{\omega}}$ gives rise to a vector field coming from a
$H_{\tilde{\omega}}$.
The following theorem answers that question:
 \begin{theo}
  The new alphabet is $\tilde{\Omega}=\Omega$; moreover for
any $\tilde{\omega}$ in $\Omega$, $\D^{(1)}_{\tilde{\omega}}$ is
Hamiltonian, and has the following expression:
\eq{\D^{(1)}_{\mot{\omega}} =  \sum_{\substack{\mot{\omega}\in\alphabet{\Omega}_r\\
\norme{\mot{\omega}}= \tilde{\omega}}} \Sam^{\mot{\omega}}
\D_{\mot{\omega}}.}
 \end{theo}
\begin{proof}
  Remember first that $\moule{\Sam}$ is alternal, so we can
still write:
  \[ \sum_{\mot{\omega}\in\alphabet{\Omega}} \Sam^{\mot{\omega}}
  \D_{\mot{\omega}} = \sum_{r\geq1} \frac{1}{r}
  \sum_{\mot{\omega}\in\alphabet{\Omega_r}} \Sam^{\mot{\omega}}
  \D_{[\mot{\omega}]}; \]
  and even
\[ \sum_{\substack{\mot{\omega}\in\alphabet{\Omega}\\
         \norme{\mot{\omega}}=\tilde{\omega}}} \Sam^{\mot{\omega}}
  \D_{\mot{\omega}} = \sum_{r\geq1} \frac{1}{r}
  \sum_{\substack{\mot{\omega}\in\alphabet{\Omega_r}\\
          \norme{\mot{\omega}}= \tilde{\omega}}} \Sam^{\mot{\omega}}
  \D_{[\mot{\omega}]}. \]
For two operators, $\D_{\omega^1}$ and $\D_{\omega^2}$, which
respectively come from two Hamiltonians
 \[H_{\omega^1}(x,y) =
\sum_{\substack{n,m\\ \omega_{nm}=\omega^1}}
a_{nm}x^ny^m \text{ and } H_{\omega^2}(x,y) = \sum_{\substack{p,q\\
\omega_{pq}=\omega^2}} a_{pq}x^py^q,\] we have actually:
 \[ [\D_{\omega^1},\D_{\omega^2}] = \poisson{
 X_{H_{\omega^1}}}{
 X_{H_{\omega^2}} } = X_{\poisson{H_{\omega^1}}{H_{\omega^2}}},\]
and \[ \poisson{H_{\omega^1}}{H_{\omega^2}} = \sum_{i=1}^\nu
\sum_{\substack{n,m,p,q\\ \omega_{nm}=\omega^1 \\ \omega_{pq} =
\omega^2}}  a_{nm}a_{pq} x^{\hat{n+p}_i}y^{\hat{m+q}_i}
(m_ip_i-q_in_i); \] this is a sum (indexed by $i$) of monomial
Hamiltonians whose each term has the same $\tilde{\omega}$:
\[\forall i, 1\leq i\leq\nu, \omega(\hat{n+p}_i,\hat{m+q}_i) =
\omega(n,m)+\omega(p,q)=\omega^1+\omega^2=\tilde{\omega}.\]
We thus can say that $\tilde{\Omega}=\left\{\norme{\mot{\omega}},
\omega\in\Omega\right\} = \Omega$. Now, remember that
$\D^{(1)}_{\tilde{\omega}}$ is defined as the Hamiltonian vector
field coming from the $\Omega$-homogeneous component of degree
$\tilde{\omega}$ of the new Hamiltonian $H^{(1)}$; we conclude
therefore that $[\D_{\omega^1},\D_{\omega^2}]$ appears in (and only
in) $\D^{(1)}_{\omega^1+\omega^2}$. Conversely, if $\tilde{\omega}$
is fixed, only the operators $\D_{[\mot{\omega}]}$ build on the words
$\mot{\omega}$ such that $\norme{\mot{\omega}}=\tilde{\omega}$ will
appear in $\D^{(1)}_{\omega}$.

Hence the result:
\eq{
\D^{(1)}_{\tilde{\omega}} = \sum_{r\geq1}
\frac{1}{r}
\sum_{\substack{\mot{\omega}\in\alphabet{\Omega}_r\\
\norme{\mot{\omega}}=\tilde{\omega}}}
\Sam^{\mot{\omega}}\D_{[\mot{\omega}]} =
\sum_{\substack{\mot{\omega}\in \Omega^*\\
\norme{\mot{\omega}}=\tilde{\omega}}} \Sam^{\mot{\omega}}
\\D_{\mot{\omega}}.
}
which concludes the proof.
\end{proof}
We may now cite the following:
\begin{theo}
The trimmed form of $X_H$ is given by 
\eq{ \Xtram=\Xlin+\sum_\bullet\moule{\Tram}\D_\bullet }
where the alphabet is $\Omega(H)$ and $\moule{\Tram}$ the mould already
defined in the previous section: for a word $\mot{\omega}$ of length $r$,
\eq{ \Tram^{\mot{\omega}} = \bigl((\moule{\Sam})^{\circ r}\bigr)^{\mot{\omega}}.}
\end{theo}
 We know that, in case of Hamiltonian vector fields, there is only one prenormal (hence normal) form.
 We have here a way to compute it; it would be interesting to compare it to other classical ways.

%% file: kolmogorov.tex
%
%
\label{section:kolmogorov}
 In this section, we use the preceding trimmed form transformations to bring a Hamiltonian vector field into Kolmogorov's normal form. The Kolmogorov's theorem
ensures the persistence of a diophantine torus of a completely integrable Hamiltonian
function under a weak perturbation. We prove this theorem in the case of a perturbation
of a special form, see \emph{infra}.

 We define the algebra $\mathcal{A}_\epsilon$ of functions $f_\epsilon(p,q):\C\times
 \C\to\C$ of the form $f_\epsilon(p,q)=\sum_{s\geq0} \epsilon^s f_s(p,q)$ where the
 $f_s$ are trigonometric polynomials in $q$, the coefficients  $f_{s,k}(p)$ of which
 being polynomials in $p$, writing $f_s(p,q)= \sum_{|k|\leq K_s}f_{s,k}(p)\expe{k\cdot q}$.
 
  We define $\mathcal{A}^1_\epsilon$ the subset of $\mathcal{A}_\epsilon$ of
  trigonometric polynomials in $q$, the coefficients of which being homogeneous
   polynomials in $p$ of degree $1$.
  
 We denote by $\der(\mathcal{A}_\epsilon)$ the set of derivations over the algebra 
$\mathcal{A}_\epsilon$, and by $\der^1(\mathcal{A}_\epsilon)$ the subset of
$\der(\mathcal{A}_\epsilon)$, of derivations $D$ of the form:
 \eq{ D= A_\epsilon(p,q)\dep_p+B_\epsilon(p,q)\dep_q,}
 with $A_\epsilon(p,q), B_\epsilon(p,q)$ in $\mathcal{A}^1_\epsilon$.
 
 Moreover, we denote by $ \der^1_r(\mathcal{A}_\epsilon)\subset 
 \der^1(\mathcal{A}_\epsilon)$ the subset of derivations
  of the form $ A_\epsilon(p,q)\dep_p + B_\epsilon(p,q)\dep_q$, $A_\epsilon$
  and $B_\epsilon$ being in $\mathcal{A}^1_\epsilon$,  of which all the coefficients of $\epsilon^s, s\leq r$ have no dependence in $q$.

  Following~\cite{giorgilli} we deal here with Hamiltonian functions from
   $\R^\nu\times\T^\nu$ to $\R$, where $\T$ is the usual torus $\R/\Z$,
   and $\freq$ a vector of $\R^\nu$ being non-resonant, of the form:
 \eqn{\label{eq:hamiltonien} H_\epsilon(p,q) = \freq\cdot p + \frac{1}{2}p^2 + \epsilon f(q).}
 The "formal" Kolmogorov theorem is then:
 \begin{theo}\label{theo:kolmo:hamilton}
 Let $H_\epsilon(p,q)$ be defined as in~\eqref{eq:hamiltonien} and
  $\freq\in\R^\nu$ being diophantine. There exists a canonical
 formal transformation $q=q^\prime + \epsilon\cdots, p=p^\prime + \epsilon\cdots$,
 which brings $H_\epsilon$ into Kolmogorov normal form:
 \eq{ H_\epsilon(p,q)=\freq\cdot p + R(p,q,\epsilon),}
 with $R(p,q,\epsilon)=O(p^2)$.
 \end{theo}
More precisely, we will prove the following theorem, denoting by 
 $X_{\epsilon}$ the Hamiltonian vector field coming from
 $H_\epsilon$, and by $X_c$ the constant vector field $\freq\cdot\dep_q$:
\begin{theo}\label{theo:kolmo:champ}
Let us suppose that $X_{\epsilon}$ has been brought into the following form:
 \eq{X_{r,\epsilon}(p,q) = \freq\cdot\dep_q + \sum_{s=1}^r \epsilon^s(a_s(p)\dep_p
      + b_s(p)\dep_q) + \sum_{s\geq r+1} \epsilon^s \biggl( \sum_{l=1}^{N_s}
         a_{s,l}(p)\expe{l\cdot q}\dep_p + \sum_{l=1}^{M_s} b_{s,l}(p) \expe{l\cdot q}
         \dep_q \biggr),}
   where $a_s(p), b_s(p), a_{s,l}(p), b_{s,l}(p)$ are in $\mathcal{A}_\epsilon^1$.
   Then, there exists a canonical transformation
 $q=q^\prime + \epsilon\cdots, p=p^\prime + \epsilon\cdots$,
 such that $X_{r,\epsilon}(p,q)=X_{r+1,\epsilon}(p^\prime,q^\prime)$.
\end{theo}
The theorem we want to prove  may be rewritten as follows:
\begin{theo}\label{theo:iteration-derivee}
If $X-X_c\in \der^1_r(\mathcal{A}_\epsilon)$, then $\Xsam-X_c \in \der^1_{r+1}(\mathcal{A}_\epsilon)$.
\end{theo}
\begin{proof}
We write 
\eq{ X = X_c + \sum_{s=1}^r \epsilon^r X_s^0(p) + \sum_{s\geq r+1} \epsilon^s 
     \biggl( X_s^0(p) + \sum_{k\in\Z^v} X_s^k(p,q) \biggr),}
     where $X_s^0(p)$ is independent of $q$ and $X_s^k(p,q)= \expe{(k\cdot q)}
     \biggl( a_{s,k}(p)\dep_p + b_{s,k}(p)\dep_q \biggr)$.
      
 We set then $B_k = \sum_{s\geq r+1} \epsilon^s X^k_s$ for $k\neq0$ and
 $B_0 = \sum_{s\geq1}\epsilon^s X_s^0$.
\begin{lemm}
For $k\in\Z^\nu$, $B_k$ is a homogeneous differential operator of degree $k$ in the 
angles $q$.  Moreover, 
\eq{ X_c B_k = \i (k\cdot\freq) B_k.}
\end{lemm}
  We can then write $X = X_c + \sum_{\mot{k}\in A^*}I^{\mot{k}} B_{\mot{k}}$,
  where $A=\Z^\nu$ is the alphabet, and $\moule{I}$ the mould already defined.
  We look for $\Theta=\sum_\bullet\moule{\Theta}\comoule{B}$, with
  $\Theta=\exp(V)$ and $V=\sum_\bullet\moule{V}\comoule{B}$.
  
  The Campbell-Baker-Hausdorff formula ensures that:
 \eqn{\label{eq:cbh}
 \Theta X \inv{\Theta} = X - [X,V] + \dotsm = X_c + \sum_{k\in\Z^\nu} B_k 
 - [X_c,V] + \hot,}
 where $\hot$ stands for higher order (in $\epsilon$) terms.
   Moreover, we set:
   \eq{ \moule{V} = \begin{cases}
                           \frac{1}{\i(k\cdot\freq)} \text{ if $\longueur{\bullet}=1$ and $k\cdot\freq
                           \neq0$;}\\
                           0 \text{ otherwise.}
                           \end{cases}}

  Then, \eq{ [X_c,V] = \sum_{k\in\Z^\nu} V^k [X_c,B_k] 
                             + \sum_{\substack{\mot{k}\in A^*\\ \longueur{\mot{k}}\geq2}}
                                     V^{\mot{k}} [X_c,B_{\mot{k}}] + \dotsm, }
    hence, \eqref{eq:cbh} rewrites:
  \eq{ X_c + 
      \underbrace{\sum_{k\in A} B_k - \sum_{k\in A}\i(\freq k )V^k B_k}_
      {(*)}
       + 
       \underbrace{\sum_{\substack{\mot{k}\in A^*\\ \longueur{\mot{k}}\geq2}} V^{\mot{k}}
         [X_c,B_{\mot{k}}] + \dotsm}_{(**)}}

The term $(**)$ is of order in $\epsilon$ at least $r+2$, therefore we do not worry
about it. The term $(*)$ rewrites:
\eq{\sum_{k\in A} \sum_{s\geq r+1} \epsilon^s X^k_s -
      \sum_{k\in A} \i(\freq\cdot k) V^k \biggl( \sum_{s\geq r+1} \epsilon^s X^k_s \biggr),}
so, if we choose $V^k= \frac{1}{\i(k\cdot\freq)}$, for $k\cdot\freq\neq0$, this latter
expression vanishes, because $X^k_{r+1}=0$ when $k\cdot\freq=0$ (as $\freq$ is non-resonant).

%
%
%
\end{proof}
 We have again our transformation $\simp$ which brings $X$ into $\Xsam$.
 After the projection lemma~\ref{lem:projection} we can write:
 \eq{ \Xsam = X_c + \sum_{r\geq1} \frac{1}{r} \sum_{\mot{k}\in A^*} \Sam^{\mot{k}} 
        B_{[\mot{k}]}.}
 Now we use the following lemma:
\begin{lemm}
\eq{ \forall D, \widetilde{D} \in \der^1(\mathcal{A}_\epsilon),\qquad
 [D,\widetilde{D}] \in \der^1(\mathcal{A}_\epsilon).}
\end{lemm}
That lemma and the projection lemma prove that $\Xsam$ is now  in
 $\der^1_{r+1}(\mathcal{A}_\epsilon)$, hence, by applying iteratively
 theorem~\ref{theo:iteration-derivee}, we are able to prove theorem~\ref{theo:kolmo:champ} therefore theorem~\ref{theo:kolmo:hamilton}.

By this way, we are able to perform a trimmed form of a Hamiltonian vector field in action-angle coordinates. It is defined by the mould $\moule{\Tram}$ exactly the same as before.

 Moreover, remark that every simplification is a canonical transformation, so at every
 step of the procedure is the vector field still Hamiltonian. That may be of great interest
 in numerical applications.

%% file: conclusion.tex

We saw in this text different powerful aspects of moulds: the ``complete calculability'' that is the universality, and the ability to be easily computed. A combinatory work (in the free Lie algebras framework) lies underneath which induces a powerful union of results both from algebra and analysis.

 The principal tool we used here was the change of graduation in the decomposition of a vector field, and we still hope to apply it to vector fields with no linear part, as E. Paul in~\cite{paul} in a future work.
 
  Moreover, the seek for normal forms has not to be limited to vector fields. We also intend to develop this kind of techniques to apply in PDEs.